\documentclass[reqno,11pt]{amsart}
\usepackage[a4paper,margin=1in]{geometry}
\usepackage[T1]{fontenc}
\usepackage{lmodern,microtype}
\usepackage{amsmath,amssymb,amsthm,amsfonts,mathrsfs,mathtools,bm}
\usepackage{epsfig, graphics, graphicx}
\usepackage{subcaption, comment, bm}
\usepackage[]{graphicx,wrapfig}
\usepackage{graphicx,booktabs,enumitem,xcolor}
\usepackage[colorlinks=true,urlcolor=magenta,linkcolor=red,citecolor=red]{hyperref}
\hypersetup{pdftitle={Planar Circular Vector Tomography with Partial Radial Data},pdfauthor={Rahul Bhardwaj and Indrani Zamindar},pdfsubject={Circular vector tomography: uniqueness, kernels, and reconstruction}}
\numberwithin{equation}{section}
\renewcommand{\baselinestretch}{1.2}
\allowdisplaybreaks[2]
\theoremstyle{plain}
\newtheorem{theorem}{Theorem}[section]
\newtheorem{proposition}[theorem]{Proposition}
\newtheorem{lemma}[theorem]{Lemma}

\theoremstyle{definition}
\newtheorem{definition}[theorem]{Definition}
\newtheorem{example}[theorem]{Example}
\theoremstyle{remark}
\newtheorem{remark}[theorem]{Remark}
\newcommand{\R}{\mathbb R}
\newcommand{\C}{\mathbb C}
\newcommand{\Z}{\mathbb Z}
\newcommand{\vx}{\boldsymbol{x}}
\newcommand{\vy}{\boldsymbol{y}}
\newcommand{\vf}{\boldsymbol{f}}
\newcommand{\vg}{\boldsymbol{g}}
\newcommand{\vn}{\boldsymbol{n}}
\newcommand{\vc}{\boldsymbol{c}}
\newcommand{\vtau}{\boldsymbol{\tau}}
\newcommand{\Lc}{\mathcal L_C}
\newcommand{\Tc}{\mathcal T_C}
\newcommand{\Rc}{\mathcal R_C}
\newcommand{\Dc}{\mathcal D_C}
\newcommand{\Vc}{\mathcal V}
\newcommand{\Yc}{\mathcal Y}
\newcommand{\DR}{\mathbb D_R}
\newcommand{\eps}{\varepsilon}
\DeclareMathOperator{\curl}{curl}
\DeclareMathOperator{\diver}{div}
\DeclareMathOperator{\supp}{supp}

\mathtoolsset{showonlyrefs=true}
\title[Circular Radon transforms for vector fields]{Circular Vector Tomography with Partial Radial Data}
\author[R. Bhardwaj]{Rahul Bhardwaj}
\address{Department of Mathematics, Indian Institute of Technology Bombay, Mumbai 400076, India}
\email{bhardwaj161067@gmail.com}
\author[I. Zamindar]{Indrani Zamindar}
\address{Johann Radon Institute for Computational and Applied Mathematics (RICAM), Altenbergerstrasse 69, 4040 Linz, Austria}
\email{indrani.zamindar@ricam.oeaw.ac.at}
\date{}
\begin{document}
\begin{abstract}
In this article, we introduce and study several generalized
circular Radon transforms acting on smooth vector fields
compactly supported in a planar annulus, with centers
restricted to the outer acquisition circle and truncated
measurement radii. The operators of interest are the
longitudinal and transverse circular Radon transforms and
their weighted analogues. We characterize
the kernels of the individual unweighted transforms. Then
we establish the reconstruction
of a vector field from the combined longitudinal and
transverse circular Radon transform data. We further show that
the pairs formed by each transform and its weighted
analogue is not sufficient for
the complete recovery of a vector field. However, we
demonstrate that full reconstruction is possible with
suitable additional data. Finally, we numerically implement
and validate the derived inversion formulas.
\end{abstract}
\maketitle

\section{Introduction}\label{sec:intro}

The circular Radon transform maps a
function to its integrals over circles.
Circular and spherical mean transforms arise in sonar and
radar imaging, thermoacoustic tomography, and photoacoustic
tomography; see, for example, \cite{AK05,AK06,FrikelQuinto15} and the
references therein. These applications have motivated
extensive research on injectivity, inversion formulas,
range characterization, and reconstruction from incomplete
data; see \cite{AK05,AK06,FPR04,Kunyansky07a}.

The reconstruction of scalar functions from circular or
spherical means has been studied using several analytical methods.
Early results include Norton's Fourier--Bessel inversion
for a circular acquisition geometry~\cite{Norton80},
Fawcett's inversion of higher-dimensional spherical
averages~\cite{Fawcett85}, and Andersson's work on the
determination of functions from spherical
averages~\cite{Andersson88}. In the
planar case, Agranovsky and Quinto~\cite{AQ} characterized
the possible sets of non-injectivity for compactly
supported functions, using microlocal methods and
properties of harmonic polynomials. Ambartsoumian and
Kuchment~\cite{AK05} subsequently studied injectivity
by techniques related to finite propagation speed and
the domain of dependence for the wave equation.

When the centers lie on the boundary of a domain
surrounding the unknown, a variety of inversion procedures
are available. Finch, and Patch~\cite{FPR04}
investigated reconstruction from spherical means with
centers on a surrounding surface, while Finch, Haltmeier,
and Rakesh~\cite{FHR07} obtained inversion formulas in
even spatial dimensions. Kunyansky developed explicit
inversion formulas~\cite{Kunyansky07a}, as well as
series-based reconstruction methods and fast
algorithms~\cite{Kunyansky07b}; later work also include
polygonal and polyhedral acquisition
surfaces~\cite{Kunyansky11}.

Alongside inversion, the characterization of data that
can arise as circular or spherical means is a fundamental
question. Ambartsoumian and Kuchment~\cite{AK06}
established a complete range description for the planar
circular transform with centers on a circle, involving
moment-type conditions and restrictions associated with
zeros of Bessel functions. Agranovsky, Kuchment, and
Quinto~\cite{AKQ07} developed range descriptions for
spherical mean transforms in higher dimensions, while
Agranovsky, Finch, and Kuchment~\cite{AFK09} investigated
the role of moment conditions in such characterizations. More recent range characterizations for
the spherical mean transform were obtained in
\cite{AAKS2023,AAKS2025,CST2026}.

A further difficulty arises when only some of the
integration circles are measured. Limited-data problems
may result from restrictions on detector positions,
accessible radii, or both. Palamodov~\cite{Palamodov00}
studied reconstruction from limited data of arc means,
and Frikel and Quinto~\cite{FrikelQuinto15} analyzed
artifacts associated with incomplete-data tomography,
including spherical mean measurements. Ambartsoumian,
Gouia-Zarrad, and Lewis~\cite{AGL} studied inversion of
the circular Radon transform on an annulus and obtained
a constructive recovery procedure when only part of the
radius range is known. Ambartsoumian and
Kunyansky~\cite{AKun14} investigated exterior and interior
circular mean problems in connection with intravascular
imaging. Ambartsoumian and Krishnan~\cite{AKrish15}
considered classes of circular and elliptical transforms
and established inversion results from radially incomplete
measurements in annular regions; related spherical-shell
reconstruction results were developed by Ambartsoumian,
Gouia-Zarrad, Krishnan, and Roy~\cite{AGKR}.
More recently, explicit inversion formulas
in odd dimensions with partial radial data were obtained
in \cite{CKT2026}.
From the computational perspective,
Haltmeier~\cite{Haltmeier16} studied angular and radial
sampling requirements for the circular Radon transform.

Vector-field tomography is motivated by imaging problems
in which directional or weighted measurements are
available. Kunyansky, McDugald, and Shearer~\cite{KMS}
studied weighted Radon transforms of vector fields in
connection with magnetoacoustoelectric tomography,
demonstrating the importance of combining unweighted
and linearly weighted observations. The extension of circular Radon transforms to vector fields is motivated by inverse problems involving directional measurements. In geometries where measurements are collected over circles, the tangential and normal components of a vector field give rise to circulation and flux measurements, respectively. Their reconstruction properties are therefore of interest in problems involving fluid flow and current-density imaging. Recent developments in photoacoustic vector tomography, such as the reconstruction of blood-flow velocity fields by Zhang et al.~\cite{Zhang2023}, further demonstrate the importance of vector-field imaging. Although these imaging methods do not directly yield circular vector integrals, they provide motivation for investigating such measurements in circular acquisition geometries. 

In this article, we study the reconstruction of smooth
vector fields compactly supported in a planar annulus
from longitudinal and transverse circular measurements,
with centers on the outer acquisition circle and a
restricted range of radii.
Building on scalar annular inversion,
we first establish constructive recovery from the joint
unweighted measurements.
We then examine the additional information carried by
first-order weighted measurements.
The longitudinal weighted pair determines the field up
to an arbitrary rotationally invariant radial component,
whereas the transverse weighted pair leaves an arbitrary
rotationally invariant tangential component undetermined.
We further show that the
undetermined component can be recovered
from an angular average of the complementary unweighted
data. Complete reconstruction is also possible using two
pointwise independent linear weights
together with the unweighted longitudinal
data. Finally, we numerically implement and validate
these formulas for
smooth vector fields compactly supported
in the annulus using MATLAB.

The article is organized as follows.
Section~\ref{sec:notations} introduces the notation and definition of generalized circular Radon transforms, and Section~\ref{sec:main}
states the main results. Section~\ref{sec:proofs}
develops the scalar reduction and the recovery from
joint unweighted measurements.
Sections~\ref{sec:longitudinal} and~\ref{sec:transverse}
consider the weighted longitudinal and transverse
transforms. Complete recovery using complementary angular
averages and two independent linear weights is
established in Sections~\ref{sec:averages}
and~\ref{sec:weights}.
Section~\ref{sec:examples} gives explicit kernel examples
and reconstruction procedures, and
Section~\ref{sec:numerical} presents numerical results.

\section{Notations and Definitions}\label{sec:notations}
We begin by introducing the notation and differential operators used throughout the article. Scalar quantities and scalar-valued functions are written in regular typeface. Points and vector fields in \(\R^2\) are written in boldface, such as \(\vx=(x_1,x_2)\) and \(\vf=(f_1,f_2)\). The Euclidean inner product is denoted by a usual dot product, and \(I_2\) is the \(2\times2\) identity matrix.
For an open set \(U\), the space \(C_c^\infty(U)\) consists of smooth functions with compact support in \(U\). A field in \(C_c^\infty(U;\R^2)\) has both components in the space \(C_c^\infty(U)\). We write \(K\Subset U\) when \(K\) is compact and contained in \(U\). 
The gradient, rotated gradient, divergence, and curl are defined by
\begin{align}
\nabla\psi&=(\partial_1\psi,\partial_2\psi),&
\nabla^\perp\psi&=(-\partial_2\psi,\partial_1\psi),\label{eq:gradients}\\
\diver\vf&=\partial_1f_1+\partial_2f_2,&
\curl\vf&=\partial_1f_2-\partial_2f_1.\label{eq:divcurl}
\end{align}
The fixed counterclockwise quarter-turn is represented by
$J=\begin{pmatrix}0&-1\\1&0\end{pmatrix}$, which satisfies
$J^T=-J$ and $J^2=-I_2$. Thus $\nabla^\perp=J\nabla$.

Let \(R>0\). The acquisition circle is
$C(R)=\{R(\cos\phi,\sin\phi):0\le\phi<2\pi\}.$
For \(\phi\in[0,2\pi)\), set
$e_\phi=(\cos\phi,\sin\phi),$ $ e_\phi^\perp=(-\sin\phi,\cos\phi),$ and $ \vc_\phi=Re_\phi.$
For \(\rho>0\), define the integration circle and disk by
\begin{align}
C(\rho,\phi)&=\{\vx:|\vx-\vc_\phi|=\rho\},\qquad \text{and}\qquad
B(\rho,\phi)=\{\vx:|\vx-\vc_\phi|<\rho\}.
\end{align}
We parametrize the circle by
\begin{equation}\label{eq:parametrization}
\vx(\alpha)=\vc_\phi+\rho(\cos\alpha,\sin\alpha),\qquad 0\le\alpha<2\pi.
\end{equation}
Since \(|\vx'(\alpha)|=\rho\), its arclength element is \(ds=\rho\,d\alpha\). The outward normal and counterclockwise tangent are
\begin{equation}\label{eq:nt}
\vn_\phi(\vx)=\frac{\vx-\vc_\phi}{|\vx-\vc_\phi|},\qquad
\vtau_\phi(\vx)=J\vn_\phi(\vx).
\end{equation}
On \eqref{eq:parametrization}, these are \((\cos\alpha,\sin\alpha)\) and \((-\sin\alpha,\cos\alpha)\), respectively. 
\begin{definition}[Scalar circular and disk transforms]
\label{def:scalar}
For $h\in C_c^\infty(\R^2)$,$\rho>0$ and $\phi\in[0,2\pi)$, the scalar circular and disk
transforms, denoted by $\Rc$ and $\Dc$, respectively,
are defined by
\begin{equation}
\Rc h(\rho,\phi)
=\int_{C(\rho,\phi)}h(\vx)\,ds_{\vx},
\qquad
\Dc h(\rho,\phi)
=\int_{B(\rho,\phi)}h(\vx)\,d\vx,
\end{equation}
\end{definition}

\begin{definition}[Longitudinal circular Radon transform]
\label{def:L}
For $\vf\in C_c^\infty(\R^2;\R^2)$, $\rho>0$ and $\phi\in[0,2\pi)$, the longitudinal
circular Radon transform, denoted by $\Lc$, is defined by
\begin{equation}\label{eq:Ldef}
\Lc\vf(\rho,\phi)
=\int_{C(\rho,\phi)}
\vf(\vx)\cdot\vtau_\phi(\vx)\,ds_{\vx},
\end{equation}
\end{definition}

\begin{definition}[Transverse circular Radon transform]
\label{def:T}
For $\vf\in C_c^\infty(\R^2;\R^2)$, $\rho>0$ and $\phi\in[0,2\pi)$, the transverse
circular Radon transform, denoted by $\Tc$, is defined by
\begin{equation}\label{eq:Tdef}
\Tc\vf(\rho,\phi)
=\int_{C(\rho,\phi)}
\vf(\vx)\cdot\vn_\phi(\vx)\,ds_{\vx},
\end{equation}
\end{definition}

\begin{definition}[First-order weighted longitudinal transform]
\label{def:L1}
For $\vf\in C_c^\infty(\R^2;\R^2)$, let
\begin{equation}\label{eq:w}
w_\phi(\vx):=(\vx-\vc_\phi)\cdot e_\phi.
\end{equation}
The first-order weighted longitudinal circular Radon
transform, denoted by $\Lc^1$, is defined by
\begin{equation}\label{eq:L1def}
\Lc^1\vf(\rho,\phi)
=\int_{C(\rho,\phi)}
w_\phi(\vx)\vf(\vx)\cdot\vtau_\phi(\vx)\,ds_{\vx},
\end{equation}
for $\rho>0$ and $\phi\in[0,2\pi)$.
\end{definition}

\begin{figure}[h!]
\centering
\includegraphics[width=.9\textwidth]{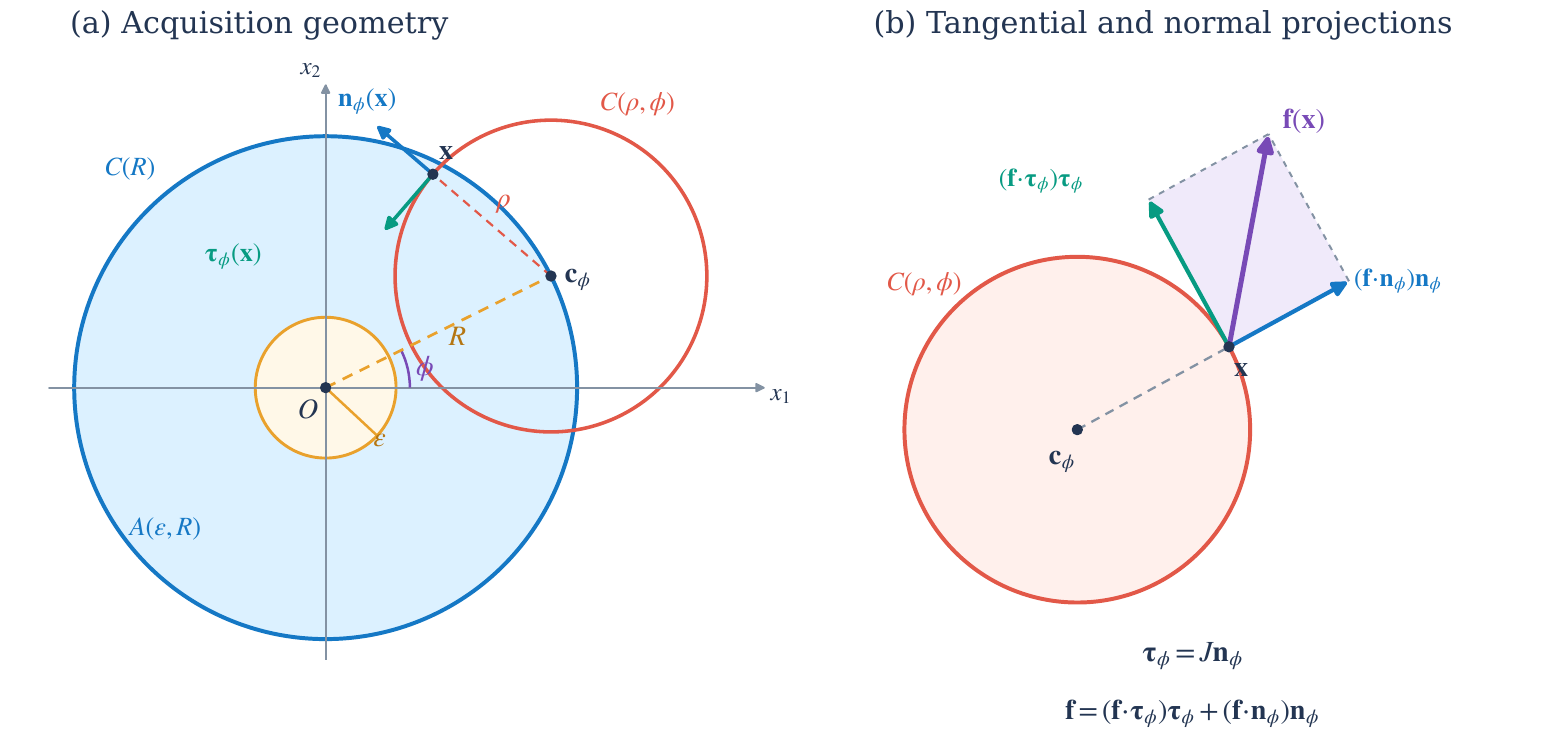}
\caption{(a) Circular acquisition geometry in the annulus
$A(\eps,R)$. (b) Tangential and normal components of $\vf$,
measured by $\Lc$ and $\Tc$, respectively.}
\label{fig:vector-circular-geometry}
\end{figure}

\begin{definition}[First-order weighted transverse transform]
\label{def:T1}
For $\vf\in C_c^\infty(\R^2;\R^2)$, $\rho>0$ and $\phi\in[0,2\pi)$ and the weight
$w_\phi$ defined in \eqref{eq:w}, the first-order weighted
transverse circular Radon transform, denoted by $\Tc^1$,
is defined by
\begin{equation}\label{eq:T1def}
\Tc^1\vf(\rho,\phi)
=\int_{C(\rho,\phi)}
w_\phi(\vx)\vf(\vx)\cdot\vn_\phi(\vx)\,ds_{\vx},
\end{equation}
\end{definition}
Fix $0<\eps<R$ and define the admissible field class and
the data domain by
\begin{equation}\label{eq:class}
\Vc_{\eps,R}:=C_c^\infty(A(\eps,R);\R^2),
\qquad
a:=R-\eps,
\qquad
\Yc:=(0,a)\times S^1.
\end{equation}
Figure~\ref{fig:vector-circular-geometry}(a) illustrates the
acquisition geometry. For $0<\rho<R-\eps$, every
$\vx\in C(\rho,\phi)$ satisfies
$|\vx|\ge|\vc_\phi|-|\vx-\vc_\phi|
=R-\rho>\eps.$
At each point of the integration circle, the outward unit
normal $\vn_\phi(\vx)$ and the counterclockwise unit tangent
$\vtau_\phi(\vx)$ form an orthonormal basis of $\R^2$.

In polar coordinates centered at the origin, let
\[
e_r(\theta)=(\cos\theta,\sin\theta),
\qquad
e_\theta(\theta)=(-\sin\theta,\cos\theta).
\]
For $\vf\in\Vc_{\eps,R}$, consider the polar decomposition
\begin{equation}\label{eq:polar-field}
\vf(r,\theta)
=f_r(r,\theta)e_r(\theta)
+f_\theta(r,\theta)e_\theta(\theta),
\end{equation}
where $f_r=\vf\cdot e_r$ and $f_\theta=\vf\cdot e_\theta$
are the radial and tangential components, respectively.
Their angular Fourier expansions are
\begin{equation}\label{eq:polar-series}
f_r(r,\theta)=\sum_{n\in\Z}a_n(r)e^{in\theta},
\qquad
f_\theta(r,\theta)=\sum_{n\in\Z}b_n(r)e^{in\theta},
\end{equation}
with coefficients
\[
a_n(r)=\frac1{2\pi}\int_0^{2\pi}
f_r(r,\theta)e^{-in\theta}\,d\theta,
\qquad
b_n(r)=\frac1{2\pi}\int_0^{2\pi}
f_\theta(r,\theta)e^{-in\theta}\,d\theta.
\]

\section{Main Results}\label{sec:main}


Throughout this section, the vector fields belong to
$\Vc_{\eps,R}$, and all measurements are restricted to
$\Yc$. Every $\vf\in\Vc_{\eps,R}$ is extended by zero
to $\R^2$, yielding a field in
$C_c^\infty(\R^2;\R^2)$, which we continue to denote
by $\vf$. Its divergence $\delta:=\diver\vf$ and curl
$\omega:=\curl\vf$ are likewise smooth and compactly
supported in $A(\eps,R)$. The kernels below are understood
with respect to this field class and data domain.

\begin{proposition}[Kernels of the unweighted transforms]
\label{thm:kernels}
The kernels of the longitudinal and transverse circular
Radon transforms are given by
\begin{align}
\ker\Lc
&=\{\vf\in\Vc_{\eps,R}:\curl\vf=0\},
\label{eq:kerL}\\
\ker\Tc
&=\{\vf\in\Vc_{\eps,R}:\diver\vf=0\}.
\label{eq:kerT}
\end{align}
\end{proposition}

\begin{theorem}[Joint unweighted recovery]
\label{thm:joint}
Every $\vf\in\Vc_{\eps,R}$ is uniquely
determined by the joint data $(\Lc\vf,\Tc\vf)$.
More precisely, scalar circular inversion determines
$\omega:=\curl\vf,
$ and $
\delta:=\diver\vf$
from $\partial_\rho\Lc\vf$ and
$\partial_\rho\Tc\vf$, respectively. 
The components of $\vf$ are then recovered
as the unique solutions of the Dirichlet problems
\begin{align}
\Delta f_1
&=\partial_1\delta-\partial_2\omega
\quad\text{in }\DR,
&
f_1|_{\partial\DR}&=0,
\label{eq:poisson1}\\
\Delta f_2
&=\partial_2\delta+\partial_1\omega
\quad\text{in }\DR,
&
f_2|_{\partial\DR}&=0.
\label{eq:poisson2}
\end{align}
\end{theorem}
\begin{theorem}[Weighted longitudinal recovery]
\label{thm:Lpair}
For every $\vf\in\Vc_{\eps,R}$, the pair
$(\Lc\vf,\Lc^1\vf)$ determines the polar Fourier
coefficients $b_n$ for all $n\in\Z$ and $a_n$ for all
$n\ne0$. The coefficient $a_0$ is the only undetermined
coefficient, and the exact joint kernel is
\begin{equation}\label{eq:Lpair-kernel}
\ker(\Lc,\Lc^1)
=\{u(r)e_r(\theta):
u\in C_c^\infty((\eps,R))\}.
\end{equation}
\end{theorem}

\begin{theorem}[Weighted transverse recovery]
\label{thm:Tpair}
For every $\vf\in\Vc_{\eps,R}$, the pair
$(\Tc\vf,\Tc^1\vf)$ determines the polar Fourier
coefficients $a_n$ for all $n\in\Z$ and $b_n$ for all
$n\ne0$. The coefficient $b_0$ is the only undetermined
coefficient, and the exact joint kernel is
\begin{equation}\label{eq:Tpair-kernel}
\ker(\Tc,\Tc^1)
=\{v(r)e_\theta(\theta):
v\in C_c^\infty((\eps,R))\}.
\end{equation}
\end{theorem}

\begin{theorem}[Recovery with complementary averaged data]
\label{thm:averages}
Every $\vf\in\Vc_{\eps,R}$ is uniquely
determined by either of the data sets
\begin{equation}\label{eq:averaged-data-sets}
\bigl(\Lc\vf,\Lc^1\vf,\overline{\Tc\vf}\bigr),
\qquad
\bigl(\Tc\vf,\Tc^1\vf,\overline{\Lc\vf}\bigr).
\end{equation}
The supplementary angular average is a function of
$\rho$. It determines $a_0$ for the longitudinal
weighted pair and $b_0$ for the transverse weighted pair.
\end{theorem}

For $b\in C^\infty(S^1;\R^2)$, define the linear weight
and the corresponding weighted longitudinal transform by
\begin{equation}\label{eq:general-weight}
\begin{aligned}
w_b(\vx,\phi)
&:=(\vx-\vc_\phi)\cdot b(\phi),\\
\mathcal L_{C,b}\vf(\rho,\phi)
&:=\int_{C(\rho,\phi)}
w_b(\vx,\phi)\vf(\vx)\cdot\vtau_\phi(\vx)\,ds_{\vx}.
\end{aligned}
\end{equation}

\begin{proposition}[Two independent linear weights]
\label{thm:two-weights}
Let $b_1,b_2\in C^\infty(S^1;\R^2)$ satisfy
\[
\det\bigl(b_1(\phi),b_2(\phi)\bigr)\ne0
\qquad\text{for every }\phi\in S^1.
\]
Then every $\vf\in\Vc_{\eps,R}$ is uniquely determined by the three data sets
$\Lc\vf,$ $
\mathcal L_{C,b_1}\vf,$ and $
\mathcal L_{C,b_2}\vf.$
\end{proposition}

\section{Proofs}\label{sec:proofs}
\subsection{Integral identities for the longitudinal and transverse transforms}
We first express both measurements as disk integrals. The precise signs depend on the orientation, which has been fixed in \eqref{eq:nt}.

\begin{proposition}\label{prop:disk-identities}
For every $\vf\in C_c^\infty(\R^2;\R^2)$, $\rho>0$,
and $\phi\in[0,2\pi)$, the longitudinal and transverse
circular Radon transforms satisfy
\begin{align}
\Lc\vf(\rho,\phi)
&=\int_{B(\rho,\phi)}\curl\vf(\vx)\,d\vx,
\label{eq:Ldisk}\\
\Tc\vf(\rho,\phi)
&=\int_{B(\rho,\phi)}\diver\vf(\vx)\,d\vx.
\label{eq:Tdisk}
\end{align}
\end{proposition}
\begin{proof}
Let \(\vx(\alpha)\) be the parametrization in \eqref{eq:parametrization}. Then
\[
\vtau_\phi(\vx(\alpha))\,ds
=(-\sin\alpha,\cos\alpha)\rho\,d\alpha
=(dx_1,dx_2).
\]
Consequently, the longitudinal circular Radon transform
\[
\Lc\vf=\int_{C(\rho,\phi)}f_1\,dx_1+f_2\,dx_2.
\]
The circle is traversed counterclockwise, so Green's theorem gives
\[
\Lc\vf=\int_{B(\rho,\phi)}(\partial_1f_2-\partial_2f_1)\,d\vx
=\int_{B(\rho,\phi)}\curl\vf\,d\vx.
\]
For the transverse data, \(\vn_\phi\) is the outward unit normal to \(B(\rho,\phi)\). The divergence theorem therefore yields
\[
\Tc\vf=\int_{\partial B(\rho,\phi)}\vf\cdot\vn_\phi\,ds
=\int_{B(\rho,\phi)}(\partial_1f_1+\partial_2f_2)\,d\vx.
\]
This is \eqref{eq:Tdisk} and completes the proof.
\end{proof}

\subsection{Connection with the scalar circular Radon transform}
\begin{lemma}[Differentiating disk integrals]
\label{lem:disk-differentiation}
For every $h\in C(\R^2)$, $\rho>0$, and
$\phi\in[0,2\pi)$, the disk integral satisfies
\begin{equation}\label{eq:disk-derivative}
\partial_\rho\int_{B(\rho,\phi)}h(\vx)\,d\vx
=\int_{C(\rho,\phi)}h(\vx)\,ds_{\vx}.
\end{equation}
\end{lemma}
\begin{proof}
Use polar coordinates centered at \(\vc_\phi\):
\[
\vx=\vc_\phi+t(\cos\alpha,\sin\alpha),\qquad
0\le t<\rho,\quad0\le\alpha<2\pi.
\]
The Jacobian is \(t\), and hence
\[
\int_{B(\rho,\phi)}h(\vx)\,d\vx
=\int_0^\rho\int_0^{2\pi}h(\vc_\phi+t e_\alpha)t\,d\alpha\,dt.
\]
Using the fundamental theorem of calculus, we have
\[
\partial_\rho\int_{B(\rho,\phi)}h\,d\vx
=\rho\int_0^{2\pi}h(\vc_\phi+\rho e_\alpha)\,d\alpha.
\]
The right side is the circular integral because \(ds=\rho\,d\alpha\).
\end{proof}

\begin{proposition}[Fundamental scalar identities]
\label{prop:fundamental}
For every $\vf\in C_c^\infty(\R^2;\R^2)$, $\rho>0$,
and $\phi\in[0,2\pi)$, the following identities hold:
\begin{align}
\partial_\rho\Lc\vf(\rho,\phi)
&=\Rc(\curl\vf)(\rho,\phi),
\label{eq:fundamental-L}\\
\partial_\rho\Tc\vf(\rho,\phi)
&=\Rc(\diver\vf)(\rho,\phi).
\label{eq:fundamental-T}
\end{align}
\end{proposition}
\begin{proof}
Apply Lemma~\ref{lem:disk-differentiation} to \eqref{eq:Ldisk} with \(h=\curl\vf\). This gives
\[
\partial_\rho\Lc\vf
=\partial_\rho\int_{B(\rho,\phi)}\curl\vf\,d\vx
=\int_{C(\rho,\phi)}\curl\vf\,ds
=\Rc(\curl\vf).
\]
Applying the lemma to \eqref{eq:Tdisk}, now with \(h=\diver\vf\), yields
\[
\partial_\rho\Tc\vf
=\partial_\rho\int_{B(\rho,\phi)}\diver\vf\,d\vx
=\int_{C(\rho,\phi)}\diver\vf\,ds
=\Rc(\diver\vf).
\]
Thus, both identities hold, completing the proof.
\end{proof}

{
\subsection{Scalar circular inversion on the annulus}
\label{subsec:scalar-inversion}

We recall the scalar annular inversion results of
Ambartsoumian, Gouia-Zarrad, and Lewis~\cite{AGL}
and record the formulas needed for the vector
reconstruction and numerical implementation.

For $h\in C_c^\infty(A(\eps,R);\C)$ and $g=\Rc h$, write
\begin{align}
h(r,\theta)
=\sum_{n\in\Z}h_n(r)e^{in\theta},
\qquad
g(\rho,\phi)
=\sum_{n\in\Z}g_n(\rho)e^{in\phi}.
\label{eq:scalar-series}
\end{align}
The Fourier coefficients are defined by
\[
h_n(r)=\frac1{2\pi}\int_0^{2\pi}
h(r,\theta)e^{-in\theta}\,d\theta,
\qquad
g_n(\rho)=\frac1{2\pi}\int_0^{2\pi}
g(\rho,\phi)e^{-in\phi}\,d\phi.
\]

{
Throughout, the first argument of $\beta(r,\rho)$ is the
distance from the origin, and the second is the measurement
radius. For $R-\rho<r<R$, the paired points in
Figure~\ref{fig:paired-intersections} have angles
$\theta_\pm=\phi\pm\beta(r,\rho)$ and are reflections
of one another across $\operatorname{span}\{e_\phi\}$.
Thus their angular offsets have opposite signs, yielding
the cosine factor in \eqref{eq:scalar-mode}.}

\begin{figure}[h!]
\centering
\includegraphics[width=.6\textwidth]
{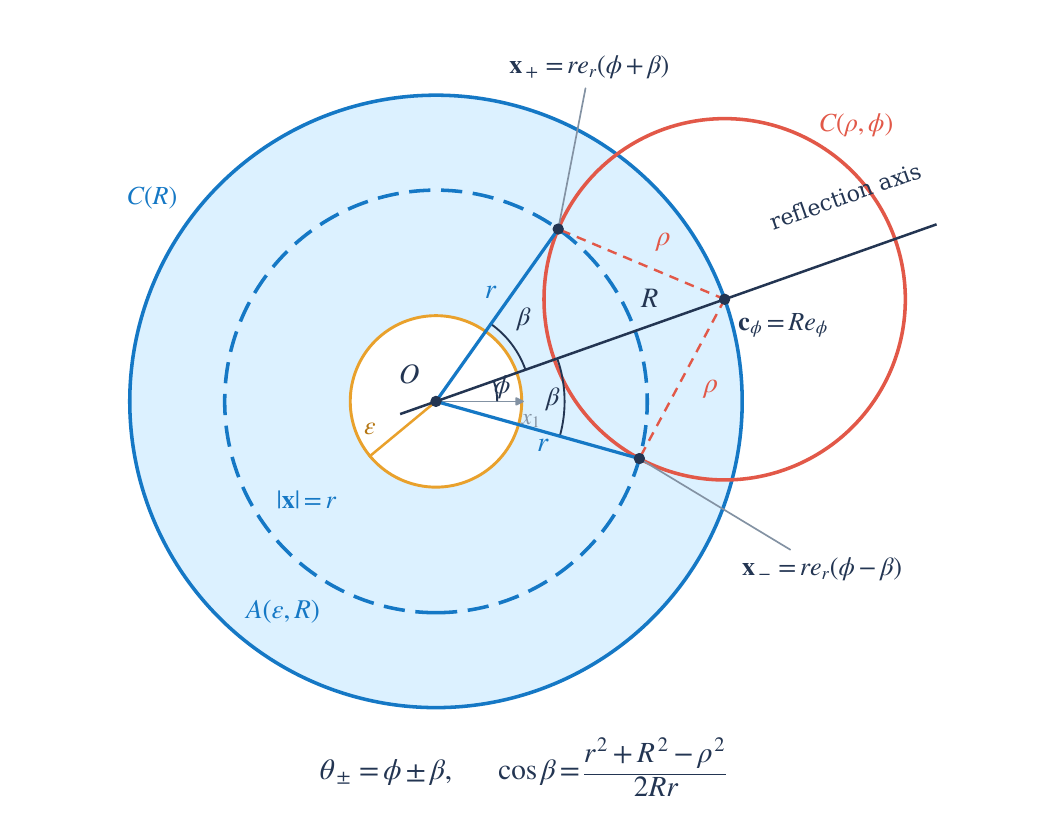}
\caption{Paired intersections $\vx_\pm=r e_r(\phi\pm{\beta(r,\rho)})$
of $|\vx|=r$ with $C(\rho,\phi)$, symmetric about
$\operatorname{span}\{e_\phi\}$. The shaded region is
the annulus $A(\eps,R)$.}
\label{fig:paired-intersections}
\end{figure}

\begin{lemma}[Scalar mode formula;
see {\cite[Section~3.2]{AGL}}]
\label{lem:scalar-mode}
For $0<\rho<R$, define
\begin{equation}\label{eq:beta}
\beta(r,\rho)
:=\arccos\left(\frac{r^2+R^2-\rho^2}{2Rr}\right),
\qquad R-\rho\le r\le R.
\end{equation}
Then
\begin{equation}\label{eq:scalar-mode}
g_n(\rho)
=\frac{2\rho}{R}\int_{R-\rho}^R
\frac{h_n(r)\cos(n\beta(r,\rho))}
{\sin\beta(r,\rho)}\,dr,
\qquad n\in\Z.
\end{equation}
On either intersection branch, parametrized by increasing
$r$, the arclength element is
\begin{equation}\label{eq:arclength-r}
ds=\frac{\rho}{R\sin\beta(r,\rho)}\,dr,
\qquad R-\rho<r<R.
\end{equation}
The integral in \eqref{eq:scalar-mode} is understood as
an improper integral at $r=R-\rho$.
\end{lemma}

\begin{lemma}[Scalar injectivity and reconstruction;
see {\cite[Theorem~1 and Corollary~2]{AGL}}]
\label{lem:scalar-inversion}
Every $h\in C_c^\infty(A(\eps,R);\C)$ is uniquely
determined by $\Rc h|_{\Yc}$. If
\[
\supp h\subset\{\vx:\alpha<|\vx|<b\},
\qquad \eps<\alpha<b<R,
\]
then each Fourier coefficient $h_n$ can be reconstructed
from $g_n$ by solving a Volterra equation of the second
kind on $[\alpha,b]$.
\end{lemma}

{We write the Fourier--Abel inversion of
\cite[Section~3.2, Corollary~2]{AGL} in our variables.
Under the support assumption of
Lemma~\ref{lem:scalar-inversion}, set $t=R-\rho$, so the
geometric angle becomes $\beta(r,R-t)$, and define}
\[
z(t,r):=\frac{r^2+R^2-(R-t)^2}{2Rr}.
\]
The identity
\begin{equation}\label{eq:sin-factorization}
4R^2r^2\sin^2\beta(r,R-t)
=(r+t)(2R+r-t)(r-t)(2R-r-t),
\end{equation}
together with
${\cos(n\beta(r,R-t))=T_{|n|}(\cos\beta(r,R-t))}$, rewrites
\eqref{eq:scalar-mode} as
\begin{equation}\label{eq:abel-first}
G_n(t):=g_n(R-t)
=\int_t^b
\frac{K_n(t,r)h_n(r)}{\sqrt{r-t}}\,dr,
\qquad \alpha\le t\le b,
\end{equation}
where
\begin{equation}\label{eq:abel-kernel}
K_n(t,r)
:=\frac{4(R-t)r\,T_{|n|}(z(t,r))}
{\sqrt{(r+t)(2R+r-t)(2R-r-t)}}.
\end{equation}
Here $T_k$ is the Chebyshev polynomial of the first kind.
For each fixed $n$, the kernel is smooth on
$\alpha\le t\le r\le b$, with diagonal value
\begin{equation}\label{eq:abel-diagonal}
d(r):=K_n(r,r)
=\sqrt{\frac{2r(R-r)}R}>0.
\end{equation}
Define
\[
u_n(r):=d(r)h_n(r),
\qquad
k_n(t,r):=\frac{K_n(t,r)}{d(r)}.
\]
For $\alpha\le s<b$, set
\begin{equation}\label{eq:second-abel}
F_n(s):=\int_s^b
\frac{G_n(t)}{\sqrt{t-s}}\,dt,
\end{equation}
and introduce
\begin{equation}\label{eq:H-kernel}
H_n(s,r)
:=\int_0^1
\frac{k_n(s+(r-s)v,r)}{\sqrt{v(1-v)}}\,dv,
\qquad \alpha\le s\le r\le b.
\end{equation}
The kernel $H_n$ is smooth, satisfies $H_n(s,s)=\pi$,
and gives
\begin{equation}\label{eq:F-H}
F_n(s)=\int_s^b H_n(s,r)u_n(r)\,dr.
\end{equation}
{Following the second-kind Volterra reduction in
\cite[Section~3.2, equation~(17)]{AGL}, with the normalization
$u_n=d\,h_n$, we obtain}
\begin{equation}\label{eq:volterra-second-kind}
u_n(s)-\int_s^b V_n(s,r)u_n(r)\,dr
=\gamma_n(s),
\end{equation}
where
\begin{equation}\label{eq:volterra-kernel-data}
V_n(s,r):=\frac1\pi\partial_sH_n(s,r),
\qquad
\gamma_n(s):=-\frac1\pi F_n'(s).
\end{equation}
{After solving \eqref{eq:volterra-second-kind},
the scalar coefficient is recovered by}
\begin{equation}\label{eq:scalar-recovered}
h_n(s)=\frac{u_n(s)}{d(s)}.
\end{equation}

For exact data, the recovered coefficients coincide
with the original coefficients, so Fourier synthesis
recovers $h$ with convergence of all derivatives on
compact sets. We apply this inversion below to the
curl, divergence, and auxiliary scalar fields.
}

\subsection{Proof of the kernels}
\begin{proof}[Proof of Proposition~\ref{thm:kernels}]
We first consider the longitudinal transform. Suppose that \(\Lc\vf(\rho,\phi)=0\) for every \((\rho,\phi)\in\Yc\). Differentiating and using \eqref{eq:fundamental-L}, we obtain
\[
0=\partial_\rho\Lc\vf(\rho,\phi)=\Rc(\curl\vf)(\rho,\phi).
\]
Differential operators do not enlarge support, so \(\curl\vf\in C_c^\infty(A(\eps,R))\). Lemma~\ref{lem:scalar-inversion} gives \(\curl\vf=0\). Thus
\[
\ker\Lc\subset\{\vf\in\Vc_{\eps,R}:\curl\vf=0\}.
\]
Conversely, if \(\curl\vf=0\), Proposition~\ref{prop:disk-identities} gives
\[
\Lc\vf(\rho,\phi)=\int_{B(\rho,\phi)}\curl\vf\,d\vx=0.
\]
This proves the reverse inclusion and establishes \eqref{eq:kerL}.

For the transverse transform, suppose that \(\Tc\vf=0\) on \(\Yc\). Equation~\eqref{eq:fundamental-T} gives \(\Rc(\diver\vf)=0\). Since \(\diver\vf\in C_c^\infty(A(\eps,R))\), scalar injectivity gives \(\diver\vf=0\). Hence
\[
\ker\Tc\subset\{\vf\in\Vc_{\eps,R}:\diver\vf=0\}.
\]
If instead \(\diver\vf=0\) is known initially, then
\[
\Tc\vf(\rho,\phi)=\int_{B(\rho,\phi)}\diver\vf\,d\vx=0.
\]
This gives the reverse inclusion and proves \eqref{eq:kerT}.
\end{proof}

\subsection{Proof of recovery from the joint unweighted data}
\begin{proof}[Proof of Proposition~\ref{thm:joint}]
By Proposition~\ref{prop:fundamental}, the data determine \(\Rc\omega\) and \(\Rc\delta\), where \(\omega=\curl\vf\) and \(\delta=\diver\vf\). Lemma~\ref{lem:scalar-inversion} reconstructs both scalar fields.  We have the following relations:
\begin{align*}
\partial_1\delta-\partial_2\omega
=\Delta f_1\qquad \text{and} \qquad 
\partial_2\delta+\partial_1\omega
=\Delta f_2.
\end{align*}
Also, \(f_1=f_2=0\) near \(\partial\DR\).

The original components supply solutions to \eqref{eq:poisson1}--\eqref{eq:poisson2}. If two solutions have the same right side and boundary trace, their difference \(u\) is harmonic with zero boundary value. Multiplying \(\Delta u=0\) by \(\overline u\), integrating by parts, and taking the real part gives
$\int_{\DR}|\nabla u|^2\,d\vx=0.$
Thus \(u\) is constant, and its zero trace makes it zero. The reconstructed components are therefore unique.

An equivalent reconstruction can be obtained using the fundamental
solution of the Laplacian on the whole plane. Let
$E(\vx)=\frac{1}{2\pi}\log|\vx|,$
so that
$\Delta E=\delta_0$
in the sense of distributions, where $\delta_0$ denotes the Dirac
distribution at the origin. Define
\[
U=E*\delta,\qquad V=E*\omega,
\qquad
\mathbf F=\nabla U+\nabla^\perp V.
\]
Since $\delta,\omega\in C_c^\infty(\R^2)$ and
$E\in L_{\mathrm{loc}}^1(\R^2)$, the potentials $U$ and $V$
are smooth. Indeed, their derivatives may be transferred to
the compactly supported smooth factors in the convolutions.
Moreover,
\[
\Delta U=(\Delta E)*\delta=\delta,
\qquad
\Delta V=(\Delta E)*\omega=\omega.
\]
These identities initially hold distributionally and hence
classically, since both sides are smooth.

With the convention
$\nabla^\perp V=(-\partial_2V,\partial_1V)$, we have
$F_1=\partial_1U-\partial_2V,$
 and 
$F_2=\partial_2U+\partial_1V.$
Consequently,
\[
\Delta F_1=\partial_1\delta-\partial_2\omega,
\qquad
\Delta F_2=\partial_2\delta+\partial_1\omega.
\]
On the other hand, the definitions
$\delta=\diver\vf$ and $\omega=\curl\vf$ give
\[
\Delta f_1=\partial_1\delta-\partial_2\omega,
\qquad
\Delta f_2=\partial_2\delta+\partial_1\omega.
\]
Thus the smooth vector field
$\mathbf W:=\mathbf F-\vf$
satisfies
$\Delta W_j=0$
in $\R^2$,\, $j=1,2$.

We next verify the behavior of $\mathbf W$ at infinity.
Choose $M>0$ such that $\supp\vf\subset B_M(0)$.
Since differentiation does not enlarge support,
$\supp\delta\cup\supp\omega\subset B_M(0).$
For $|\vx|>2M$ and $\vy\in B_M(0)$, we have
$|\vx-\vy|\ge |\vx|-|\vy|\ge \frac{|\vx|}{2}.$
Using
$\nabla E(\vx)=\frac{\vx}{2\pi|\vx|^2},$
we obtain
\begin{align*}
|\nabla U(\vx)|
&\le
\int_{B_M(0)}
\frac{|\delta(\vy)|}{2\pi|\vx-\vy|}\,d\vy
\le
\frac{\|\delta\|_{L^1(\R^2)}}{\pi|\vx|},\\
|\nabla^\perp V(\vx)|
&\le
\int_{B_M(0)}
\frac{|\omega(\vy)|}{2\pi|\vx-\vy|}\,d\vy
\le
\frac{\|\omega\|_{L^1(\R^2)}}{\pi|\vx|}.
\end{align*}
Therefore,
\[
|\mathbf F(\vx)|
\le
\frac{
\|\delta\|_{L^1(\R^2)}
+\|\omega\|_{L^1(\R^2)}
}{\pi|\vx|}
\longrightarrow 0
\qquad\text{as }|\vx|\to\infty.
\]
Because $\vf$ is compactly supported, it follows that
$\mathbf W(\vx)\to0$ as $|\vx|\to\infty$.

Each component $W_j$ is therefore a bounded harmonic function
on $\R^2$. By the harmonic Liouville theorem, $W_j$ is constant,
and its vanishing limit at infinity makes this constant zero.
Hence $\mathbf F=\vf$, proving the whole-plane reconstruction
formula
\begin{equation}\label{eq:whole-plane-recovery}
\vf
=
\nabla(E*\delta)+\nabla^\perp(E*\omega),
\qquad
E(\vx)=\frac{1}{2\pi}\log|\vx|.
\end{equation}
\end{proof}
\section{Recovery of \texorpdfstring{$\vf$}{the field} from \texorpdfstring{$\Lc\vf$ and $\Lc^1\vf$}{longitudinal weighted data}}\label{sec:longitudinal}
In this section we prove Theorem~\ref{thm:Lpair} constructively. All formulas are written for the data convention \((\rho,\phi)\) used in the definitions.
\subsection{The curl equations for the polar Fourier coefficients}
The first identity in Proposition~\ref{prop:fundamental} and scalar inversion determine
$\omega(\vx):=\curl\vf(\vx)$
from \(\Lc\vf\). Let \(\omega(r,\theta)=\sum_{n\in\Z}\omega_n(r)e^{in\theta}\). We first obtain the precise relation between the known quantity and the unknown polar coefficients.

In polar coordinates, the Cartesian derivatives are
\begin{equation}\label{eq:cartesian-polar-derivatives}
\partial_{x_1}=\cos\theta\,\partial_r-\frac{\sin\theta}{r}\partial_\theta,
\qquad
\partial_{x_2}=\sin\theta\,\partial_r+\frac{\cos\theta}{r}\partial_\theta.
\end{equation}
The Cartesian components of the vector field $\vf$
in \eqref{eq:polar-field} are
$f_1=f_r\cos\theta-f_\theta\sin\theta,
$ and $
f_2=f_r\sin\theta+f_\theta\cos\theta.$
Substitution into \(\partial_1f_2-\partial_2f_1\), followed by collecting the radial and angular derivatives, gives
\begin{align}
\omega
&=\partial_rf_\theta+\frac{f_\theta}{r}-\frac1r\partial_\theta f_r\nonumber\\
&=\frac1r\left[\partial_r(rf_\theta)-\partial_\theta f_r\right].
\label{eq:polar-curl}
\end{align}
Substituting the expansion for $\omega(r,\theta)$ together with
\eqref{eq:polar-series} into \eqref{eq:polar-curl}
and comparing the coefficients of $e^{in\theta}$,
we obtain
\begin{equation}\label{eq:curl-mode}
(rb_n)'(r)-in\,a_n(r)=r\omega_n(r),\qquad n\in\Z.
\end{equation}
For \(n\ne0\), this determines \(a_n\) once \(b_n\) is known:
\begin{equation}\label{eq:a-from-b}
a_n(r)=\frac{(rb_n)'(r)-r\omega_n(r)}{in}.
\end{equation}
It remains to recover \(b_n\) for each nonzero mode.

\subsection{Removing the known weighted curl term}
 Using \eqref{eq:w}, for fixed \(\phi\), \(\partial_1w_\phi=\cos\phi\) and \(\partial_2w_\phi=\sin\phi\). Hence
\begin{align}
\curl(w_\phi\vf)
&=\partial_1(w_\phi f_2)-\partial_2(w_\phi f_1)\nonumber\\
&=w_\phi(\partial_1f_2-\partial_2f_1)
 +(\partial_1w_\phi)f_2-(\partial_2w_\phi)f_1\nonumber\\
&=w_\phi\omega+f_2\cos\phi-f_1\sin\phi\nonumber\\
&=w_\phi\omega+\vf\cdot e_\phi^\perp.
\label{eq:weighted-curl-product}
\end{align}
Green's theorem applied to \(w_\phi\vf\) in the definition \eqref{eq:L1def} now gives
\begin{equation}\label{eq:L1-disk}
\Lc^1\vf(\rho,\phi)
=\int_{B(\rho,\phi)}w_\phi(\vx)\omega(\vx)\,d\vx
 +\int_{B(\rho,\phi)}\vf(\vx)\cdot e_\phi^\perp\,d\vx.
\end{equation}
The first term is known, since the field \(\omega\) has already been reconstructed. Define the known residual
\begin{equation}\label{eq:QL}
Q_L(\rho,\phi):=\Lc^1\vf(\rho,\phi)
-\int_{B(\rho,\phi)}w_\phi(\vx)\omega(\vx)\,d\vx.
\end{equation}
Then
\begin{equation}\label{eq:QL-disk}
Q_L(\rho,\phi)=\int_{B(\rho,\phi)}\vf(\vx)\cdot e_\phi^\perp\,d\vx.
\end{equation}
For each fixed \(\phi\), the vector \(e_\phi^\perp\) is independent of \(\rho\) and \(\vx\). Lemma~\ref{lem:disk-differentiation} therefore gives
\begin{equation}\label{eq:PL}
P_L(\rho,\phi):=\partial_\rho Q_L(\rho,\phi)
=\int_{C(\rho,\phi)}\vf(\vx)\cdot e_\phi^\perp\,ds_{\vx}.
\end{equation}


\begin{lemma}\label{lem:longitudinal-mode} Write \(P_L(\rho,\phi)=\sum_{n\in\Z}P_{L,n}(\rho)e^{in\phi}\).
For every \(n\in\Z\), the Fourier
coefficient $P_{L,n}$ satisfies
\begin{equation}\label{eq:PLmode}
P_{L,n}(\rho)=\frac{2\rho}{R}\int_{R-\rho}^R
\left[ia_n(r)\sin(n{\beta(r,\rho)})
+b_n(r)\frac{\cos{\beta(r,\rho)}\cos(n{\beta(r,\rho)})}{\sin{\beta(r,\rho)}}\right]dr,
\end{equation}
{Here $\beta(r,\rho)$ is defined by \eqref{eq:beta}.}
\end{lemma}
\begin{proof}
The projections of the polar basis are
$e_r(\theta)\cdot e_\phi^\perp=\sin(\theta-\phi),
$ and $
e_\theta(\theta)\cdot e_\phi^\perp=\cos(\theta-\phi).$
Thus, we have
\[
\vf(r,\theta)\cdot e_\phi^\perp
=f_r(r,\theta)\sin(\theta-\phi)
 +f_\theta(r,\theta)\cos(\theta-\phi).
\]
For a fixed interior radius \(r\), the two intersection points consequently have angles \(\theta_\pm=\phi\pm{\beta(r,\rho)}\). They are reflections of each other across the line joining the origin to \(c_\phi\).
This reflection is depicted in
Figure~\ref{fig:paired-intersections}, where the two points
$\vx_+$ and $\vx_-$ correspond to the angles
$\theta_+=\phi+{\beta(r,\rho)}$ and $\theta_-=\phi-{\beta(r,\rho)}$, respectively. After adding the two contributions of the \(n\)-th radial mode, we get
\begin{align*}
&a_n(r)e^{in(\phi+{\beta(r,\rho)})}\sin{\beta(r,\rho)}
+a_n(r)e^{in(\phi-{\beta(r,\rho)})}\sin(-{\beta(r,\rho)})\\
&\qquad=a_n(r)e^{in\phi}\sin{\beta(r,\rho)}
\bigl(e^{in{\beta(r,\rho)}}-e^{-in{\beta(r,\rho)}}\bigr)\\
&=2i\,a_n(r)e^{in\phi}\sin{\beta(r,\rho)}\sin(n{\beta(r,\rho)}).
\end{align*}
The contributions of the \(n\)-th tangential mode are
\begin{align*}
&b_n(r)e^{in(\phi+{\beta(r,\rho)})}\cos{\beta(r,\rho)}
+b_n(r)e^{in(\phi-{\beta(r,\rho)})}\cos(-{\beta(r,\rho)})\\
&\qquad=b_n(r)e^{in\phi}\cos{\beta(r,\rho)}
\bigl(e^{in{\beta(r,\rho)}}+e^{-in{\beta(r,\rho)}}\bigr)\\
&=2b_n(r)e^{in\phi}\cos{\beta(r,\rho)}\cos(n{\beta(r,\rho)}).
\end{align*}
Multiplying their sum by the single-branch arclength element \eqref{eq:arclength-r}, we obtain
\[
\frac{2\rho}{R}e^{in\phi}
\left[ia_n(r)\sin(n{\beta(r,\rho)})
+b_n(r)\frac{\cos{\beta(r,\rho)}\cos(n{\beta(r,\rho)})}{\sin{\beta(r,\rho)}}\right]dr.
\]
Integrating over \(R-\rho<r<R\) and comparing Fourier coefficients proves the formula.
\end{proof}

\subsection{Integration by parts and cancellation of the geometric factor}
Fix \(n\ne0\), and let
$y_{L,n}(r):=rb_n(r).$
From \eqref{eq:curl-mode}, we have
\[
ia_n(r)=\frac{y_{L,n}'(r)-r\omega_n(r)}n.
\]
Substituting this into \eqref{eq:PLmode} yields
\begin{align}
\frac{R}{2\rho}P_{L,n}(\rho)
={}&\frac1n\int_{R-\rho}^R y_{L,n}'(r)\sin(n{\beta(r,\rho)})\,dr-\frac1n\int_{R-\rho}^R r\omega_n(r)\sin(n{\beta(r,\rho)})\,dr\nonumber\\
&\quad+\int_{R-\rho}^R\frac{y_{L,n}(r)}r
\frac{\cos{\beta(r,\rho)}\cos(n{\beta(r,\rho)})}{\sin{\beta(r,\rho)}}\,dr.
\label{eq:PLsub}
\end{align}
We next justify the integration by parts at the lower endpoint,
where $\partial_r{\beta(r,\rho)}$ may be singular. Fix $0<\rho<R$ and
$n\in\mathbb Z\setminus\{0\}$, and put $t=R-\rho>0$.
Recall that
\[
\cos\beta(r,\rho)
=\frac{r^2+R^2-\rho^2}{2Rr},
\qquad t\leq r\leq R.
\]
In particular, $\beta(t,\rho)=0$. By
\eqref{eq:sin-factorization}, we have
\[
\sin^2\beta(r,\rho)
=
\frac{(r-t)(r+t)\bigl((R+\rho)^2-r^2\bigr)}
     {4R^2r^2}.
\]
Thus, for $r\geq t$,
\[
\sin\beta(r,\rho)=\sqrt{r-t}\,q_\rho(r),
\qquad
q_\rho(r):=
\frac{\sqrt{(r+t)\bigl((R+\rho)^2-r^2\bigr)}}{2Rr}.
\]
The function $q_\rho$ is smooth and strictly positive on $[t,R]$,
and
$q_\rho(t)=\sqrt{\frac{2\rho}{Rt}}>0.$
Consequently, there exist $h>0$ and constants $c,C>0$ such that
\[
c\sqrt{r-t}\leq \sin\beta(r,\rho)
\leq C\sqrt{r-t},
\qquad t<r<t+h.
\]
Since $\beta(r,\rho)\to0$ as $r\downarrow t$ and
$\sin u/u\to1$ as $u\to0$, it also follows that
\[
\beta(r,\rho)=O\bigl(\sqrt{r-t}\bigr).
\]
Moreover, differentiating the identity for $\cos{\beta(r,\rho)}$ gives
\[
-\sin{\beta(r,\rho)}\,\partial_r{\beta(r,\rho)}
=\frac{r^2-R^2+\rho^2}{2Rr^2},
\]
and hence
\[
\partial_r\beta(r,\rho)
=\frac{R^2-\rho^2-r^2}{2Rr^2\sin\beta(r,\rho)}
=O\bigl((r-t)^{-1/2}\bigr).
\]
Here and below, the endpoint estimates are understood as
$r\downarrow t$, with $\rho$ and $n$ fixed. In particular,
\[
|\sin(n{\beta(r,\rho)})|
\leq |n|{\beta(r,\rho)}
=O\bigl(\sqrt{r-t}\bigr).
\]

Set $y(r)=y_{L,n}(r)=r b_n(r)$. The smoothness of the field
implies that $y$ and $y'$ are bounded on $[t,R]$. Furthermore,
compact support of the field inside $D_R$ implies that $y$
vanishes in a neighbourhood of $R$.
For $0<\eta<R-t$, ordinary integration by parts on
$[t+\eta,R]$ yields
\begin{align*}
\frac1n\int_{t+\eta}^R y'(r)\sin(n{\beta(r,\rho)})\,dr
&=
\left[\frac{y(r)\sin(n{\beta(r,\rho)})}n\right]_{t+\eta}^{R}
-\int_{t+\eta}^R
y(r)\cos(n{\beta(r,\rho)})\,\partial_r{\beta(r,\rho)}\,dr,
\end{align*}
where we used
$\partial_r\sin(n{\beta(r,\rho)})
=n\cos(n{\beta(r,\rho)})\,\partial_r{\beta(r,\rho)}.$

The boundary contribution at $R$ is zero. At the lower
endpoint, the preceding estimates give
\[
\left|
\frac{y(t+\eta)\sin(n\beta(t+\eta,\rho))}{n}
\right|
\leq \|y\|_{L^\infty([t,R])}\,
      \beta(t+\eta,\rho)
\leq C\sqrt{\eta},
\]
so this boundary contribution tends to zero as
$\eta\downarrow0$.

Finally, both integrals admit absolutely convergent limits
at $t$. Indeed, near $t$,
\[
\left|\frac1n y'(r)\sin(n{\beta(r,\rho)})\right|
\leq C\sqrt{r-t},
\qquad
\left|y(r)\cos(n{\beta(r,\rho)})\,\partial_r{\beta(r,\rho)}\right|
\leq C(r-t)^{-1/2}.
\]
More explicitly, for sufficiently small $\eta>0$,
\[
\left|
\frac1n\int_t^{t+\eta}y'(r)\sin(n{\beta(r,\rho)})\,dr
\right|
\leq C\eta^{3/2},
\qquad
\left|
\int_t^{t+\eta}
y(r)\cos(n{\beta(r,\rho)})\,\partial_r{\beta(r,\rho)}\,dr
\right|
\leq C\sqrt{\eta}.
\]
Letting $\eta\downarrow0$ in the truncated integration-by-parts
identity therefore gives
\begin{equation}\label{eq:L-ibp}
\frac1n\int_t^R
y_{L,n}'(r)\sin(n\beta(r,\rho))\,dr
=
-\int_t^R
y_{L,n}(r)\cos(n\beta(r,\rho))\,
\partial_r\beta(r,\rho)\,dr.
\end{equation}
Combining \eqref{eq:PLsub} and \eqref{eq:L-ibp}, we have
\begin{align}
&\frac{R}{2\rho}P_{L,n}(\rho)
+\frac1n\int_{R-\rho}^R r\omega_n(r)\sin(n{\beta(r,\rho)})\,dr
\nonumber\\
&\qquad=\int_{R-\rho}^R y_{L,n}(r)\cos(n{\beta(r,\rho)})
\left[-\partial_r{\beta(r,\rho)}+\frac{\cos{\beta(r,\rho)}}{r\sin{\beta(r,\rho)}}\right]dr.
\label{eq:L-before-cancellation}
\end{align}
{It remains to simplify the bracket. Differentiating
\[
\cos\beta(r,\rho)=\frac{r^2+R^2-\rho^2}{2Rr}
\]
gives
\[
-\sin\beta(r,\rho)\,\partial_r\beta(r,\rho)
=\frac{r^2-R^2+\rho^2}{2Rr^2}.
\]}
Hence
\[
-\partial_r{\beta(r,\rho)}=\frac{r^2-R^2+\rho^2}{2Rr^2\sin{\beta(r,\rho)}},
\qquad
\frac{\cos{\beta(r,\rho)}}{r\sin{\beta(r,\rho)}}
=\frac{r^2+R^2-\rho^2}{2Rr^2\sin{\beta(r,\rho)}}.
\]
Adding these two fractions yields the cancellation
\begin{equation}\label{eq:geometric-cancellation}
-\partial_r{\beta(r,\rho)}+\frac{\cos{\beta(r,\rho)}}{r\sin{\beta(r,\rho)}}
=\frac1{R\sin{\beta(r,\rho)}}.
\end{equation}
Define the known corrected coefficient
\begin{equation}\label{eq:GL}
G_{L,n}(\rho):=\frac{R}{2\rho}P_{L,n}(\rho)
+\frac1n\int_{R-\rho}^R r\omega_n(r)\sin(n\beta(r,\rho))\,dr.
\end{equation}
Using equations~\eqref{eq:L-before-cancellation} and \eqref{eq:geometric-cancellation}, we have
\begin{equation}\label{eq:L-reduced}
G_{L,n}(\rho)=\frac1R\int_{R-\rho}^R
\frac{y_{L,n}(r)\cos(n\beta(r,\rho))}{\sin\beta(r,\rho)}\,dr.
\end{equation}
Introduce the scalar field
\[
h(r,\theta):=\sum_{k\ne0}y_{L,k}(r)e^{ik\theta}
=r\bigl(f_\theta(r,\theta)-b_0(r)\bigr).
\]
Since the original field is smooth and compactly supported in
$A(\eps,R)$, the same holds for $h$. Its $n$-th angular Fourier
coefficient is $h_n=y_{L,n}$ for every $n\ne0$. Therefore,
by \eqref{eq:scalar-mode},
\begin{align}
\bigl(\Rc h\bigr)_n(\rho)
&=\frac{2\rho}{R}\int_{R-\rho}^{R}
\frac{y_{L,n}(r)\cos(n\beta(r,\rho))}
{\sin\beta(r,\rho)}\,dr \nonumber\\
&=2\rho\,G_{L,n}(\rho),\qquad n\ne0.
\label{eq:L-scalar-normalization}
\end{align}
{Applying Lemma~\ref{lem:scalar-inversion}, using
\eqref{eq:volterra-second-kind} and \eqref{eq:scalar-recovered}
with data coefficient $2\rho G_{L,n}$, determines $y_{L,n}$
for each $n\ne0$.}
  Consequently, all nonzero polar Fourier coefficients are determined
by the measured data through
\begin{equation}\label{eq:L-reconstruction-coefficients}
b_n(r)=\frac{y_{L,n}(r)}r,
\qquad
a_n(r)=\frac{y_{L,n}'(r)-r\omega_n(r)}{in},
\qquad n\ne0.
\end{equation}

\subsection{The zeroth mode and the complete joint kernel}
For \(n=0\), equation~\eqref{eq:curl-mode} gives
$(rb_0)'(r)=r\omega_0(r).$
The field vanishes near the inner boundary, so \(b_0\) vanishes there and the integration constant is zero. Therefore
\begin{equation}\label{eq:L-b0}
b_0(r)=\frac1r\int_\eps^r s\omega_0(s)\,ds.
\end{equation}
The coefficient \(a_0\) does not occur in the equation \eqref{eq:curl-mode}. In \eqref{eq:PLmode}, its coefficient is \(\sin(0{\beta(r,\rho)})=0\). Thus no step of this reconstruction determines \(a_0\). The reconstructed field is
\begin{equation}\label{eq:L-synthesis}
\vf(r,\theta)=a_0(r)e_r(\theta)
+\sum_{n\ne0}a_n(r)e^{in\theta}e_r(\theta)
+\sum_{n\in\Z}b_n(r)e^{in\theta}e_\theta(\theta),
\end{equation}
where every term except the first is determined by the data.

\begin{proof}[Proof of Theorem~\ref{thm:Lpair}]
The formulas above provide the claimed constructive recovery. To prove that they give the exact kernel, suppose \(\Lc\vf=\Lc^1\vf=0\). Then \(\omega=0\), \(Q_L=0\), and \(P_L=0\). Every \(G_{L,n}\), \(n\ne0\), is zero. Scalar injectivity applied to \eqref{eq:L-scalar-normalization} gives \(y_{L,n}=0\). Formula~\eqref{eq:L-reconstruction-coefficients} gives \(a_n=b_n=0\) for \(n\ne0\), and \eqref{eq:L-b0} gives \(b_0=0\). Consequently, \(\vf=a_0(r)e_r\), with \(a_0\in C_c^\infty((\eps,R))\).

For the reverse inclusion, take \(\vf=u(r)e_r\). Then the polar-coordinate formula
\eqref{eq:polar-curl} gives $\curl\vf=0$. Hence
\eqref{eq:Ldisk} implies $\Lc\vf=0$. Equation~\eqref{eq:L1-disk} reduces to
\[
\Lc^1\vf(\rho,\phi)=\int_{B(\rho,\phi)}u(|\vx|)e_r(\theta)\cdot e_\phi^\perp\,d\vx.
\]
{Let $S_\phi$ denote reflection across the line through the origin
and $\vc_\phi$, namely
$S_\phi\vx=2(\vx\cdot e_\phi)e_\phi-\vx.$
This reflection satisfies
$S_\phi\vc_\phi=\vc_\phi,\,
S_\phi e_\phi^\perp=-e_\phi^\perp,\,
|S_\phi\vx|=|\vx|,
\, |\det S_\phi|=1.$
Moreover, we have
$|S_\phi\vx-\vc_\phi|
=|S_\phi(\vx-\vc_\phi)|
=|\vx-\vc_\phi|,$
so $S_\phi$ maps $B(\rho,\phi)$ onto itself.

For $\vy\ne0$, the radial unit vector satisfies
$e_r(S_\phi\vy)
={S_\phi\vy}/{|\vy|}
=S_\phi e_r(\vy).$
Since $S_\phi^T=S_\phi$, it follows that
\[
e_r(S_\phi\vy)\cdot e_\phi^\perp
=e_r(\vy)\cdot S_\phi e_\phi^\perp
=-e_r(\vy)\cdot e_\phi^\perp.
\]
The integrand is understood to be zero near the origin,
where $u$ vanishes. Changing variables $\vx=S_\phi\vy$
therefore gives
\begin{align*}
\Lc^1\vf(\rho,\phi)
&=\int_{B(\rho,\phi)}
u(|S_\phi\vy|)
e_r(S_\phi\vy)\cdot e_\phi^\perp\,d\vy\\
&=-\int_{B(\rho,\phi)}
u(|\vy|)e_r(\vy)\cdot e_\phi^\perp\,d\vy\\
&=-\Lc^1\vf(\rho,\phi).
\end{align*}
Consequently, $\Lc^1\vf(\rho,\phi)=0$ for every measured
$(\rho,\phi)$. Together with $\Lc\vf=0$, this proves
\[
\{u(r)e_r(\theta):u\in C_c^\infty((\eps,R))\}
\subseteq\ker(\Lc,\Lc^1),
\]
completing the reverse inclusion. This completes the proof of Theorem~\ref{thm:Lpair}}
\end{proof}


\section{Transverse weighted recovery by a fixed rotation}\label{sec:transverse}
The longitudinal proof contains the endpoint integration by parts and the scalar normalization needed for both weighted pairs. We now transfer that proof by an explicit rotation and record every quantity needed to reconstruct the transverse data, including its signs and missing mode.

\subsection{Rotation identities and the divergence equations}
Let $\vg=J\vf$. Using $J^TJ=I_2$, $\vtau_\phi=J\vn_\phi$,
and $J^T=-J$, we obtain
\[
(J\vf)\cdot\vtau_\phi
=(J\vf)\cdot(J\vn_\phi)
=\vf\cdot\vn_\phi,
\qquad
(J\vf)\cdot\vn_\phi
=\vf\cdot J^T\vn_\phi
=-\vf\cdot\vtau_\phi.
\]
Therefore, the definitions of the longitudinal and transverse
circular Radon transforms give
\begin{align}
\Lc\vg&=\Tc\vf,
&
\Lc^1\vg&=\Tc^1\vf,
&
\Tc\vg&=-\Lc\vf.
\label{eq:transverse-rotation}
\end{align}
The weighted identity follows by multiplying the first
pointwise identity by the common scalar weight $w_\phi$
and integrating over $C(\rho,\phi)$.

Writing $\vf=(f_1,f_2)$, we have $\vg=J\vf=(-f_2,f_1)$.
Consequently,
\[
\curl\vg
=\partial_1f_1+\partial_2f_2
=\diver\vf=:\delta.
\]
Moreover, since $Je_r=e_\theta$ and $Je_\theta=-e_r$,
\[
\vg
=J(f_r e_r+f_\theta e_\theta)
=-f_\theta e_r+f_r e_\theta.
\]
Thus the polar Fourier coefficients of $\vg$ satisfy
\begin{equation}\label{eq:transverse-mode-transfer}
\widetilde a_n=-b_n,
\qquad
\widetilde b_n=a_n,
\qquad
\widetilde\omega_n=\delta_n.
\end{equation}
Applying \eqref{eq:curl-mode} to $\vg$ gives
\[
(r\widetilde b_n)'(r)-in\,\widetilde a_n(r)
=r\widetilde\omega_n(r).
\]
Substituting \eqref{eq:transverse-mode-transfer}, we obtain
\begin{equation}\label{eq:div-mode}
(ra_n)'(r)+in\,b_n(r)=r\delta_n(r),
\qquad n\in\Z.
\end{equation}
Equivalently, the polar divergence identity is
\begin{equation}\label{eq:polar-divergence}
\delta=\frac1r\bigl[\partial_r(rf_r)+\partial_\theta f_\theta\bigr].
\end{equation}
The class \(\Vc_{\eps,R}\) is preserved by \(J\), so this substitution preserves every support and smoothness hypothesis of the longitudinal proof.
\subsection{Known residuals and nonzero-mode reconstruction}

First reconstruct $\delta=\diver\vf$ from
$\partial_\rho\Tc\vf=\Rc\delta$ using the scalar
Fourier--Abel inverse. Applying \eqref{eq:L1-disk} to
$\vg=J\vf$ and using
\[
\Lc^1\vg=\Tc^1\vf,
\qquad
\curl\vg=\delta,
\qquad
\vg\cdot e_\phi^\perp
=(J\vf)\cdot(Je_\phi)=\vf\cdot e_\phi,
\]
we obtain
\begin{equation}\label{eq:T1-disk}
\Tc^1\vf(\rho,\phi)
=\int_{B(\rho,\phi)}w_\phi(\vx)\delta(\vx)\,d\vx
+\int_{B(\rho,\phi)}\vf(\vx)\cdot e_\phi\,d\vx.
\end{equation}
Consequently, the transverse weighted pair determines
\begin{align}
Q_T(\rho,\phi)
&:=\Tc^1\vf(\rho,\phi)
-\int_{B(\rho,\phi)}w_\phi(\vx)\delta(\vx)\,d\vx,
\label{eq:QT}\\
P_T(\rho,\phi)&:=\partial_\rho Q_T(\rho,\phi).
\label{eq:PTdef}
\end{align}
By the disk differentiation lemma,
\begin{equation}\label{eq:PT}
Q_T=\int_{B(\rho,\phi)}\vf\cdot e_\phi\,d\vx,
\qquad
P_T=\int_{C(\rho,\phi)}\vf\cdot e_\phi\,ds.
\end{equation}
These residuals coincide with the longitudinal residuals
of the rotated field $\vg=J\vf$. Indicating the field
argument explicitly, we have
\[
Q_T[\vf]=Q_L[\vg],
\qquad
P_T[\vf]=P_L[\vg].
\]

Write
\[
P_T(\rho,\phi)=\sum_{n\in\Z}P_{T,n}(\rho)e^{in\phi},
\qquad
\delta(r,\theta)=\sum_{n\in\Z}\delta_n(r)e^{in\theta}.
\]
Equation~\eqref{eq:GL}, applied to $\vg$ with
$\widetilde\omega_n=\delta_n$, gives the known corrected
coefficient
\begin{equation}\label{eq:GT}
G_{T,n}(\rho)
:=\frac{R}{2\rho}P_{T,n}(\rho)
+\frac1n\int_{R-\rho}^R
r\delta_n(r)\sin(n\beta(r,\rho))\,dr,
\qquad n\ne0.
\end{equation}
The coefficient reconstructed in the longitudinal
argument for $\vg$ is
\begin{equation}\label{eq:yT}
y_{T,n}(r):=r\widetilde b_n(r)=ra_n(r),
\qquad n\ne0.
\end{equation}
Since $y_{T,n}$ is smooth and compactly supported in
$(\eps,R)$, the endpoint estimates and integration-by-parts
identity \eqref{eq:L-ibp} apply to this coefficient.
The geometric cancellation
\eqref{eq:geometric-cancellation} therefore yields
\begin{equation}\label{eq:T-reduced}
G_{T,n}(\rho)
=\frac1R\int_{R-\rho}^R
\frac{y_{T,n}(r)\cos(n\beta(r,\rho))}
{\sin\beta(r,\rho)}\,dr,
\qquad n\ne0.
\end{equation}
Introduce the scalar field
\[
h_T(r,\theta)
:=\sum_{k\ne0}y_{T,k}(r)e^{ik\theta}
=r\bigl(f_r(r,\theta)-a_0(r)\bigr).
\]
This field is smooth and compactly supported in $A(\eps,R)$.
By \eqref{eq:scalar-mode} and \eqref{eq:T-reduced},
\[
(\Rc h_T)_n(\rho)=2\rho G_{T,n}(\rho),
\qquad n\ne0.
\]
{Applying Lemma~\ref{lem:scalar-inversion}, using
\eqref{eq:volterra-second-kind} and \eqref{eq:scalar-recovered}
with data coefficient $2\rho G_{T,n}$, determines $y_{T,n}$
for every $n\ne0$.} Using \eqref{eq:div-mode}, we therefore recover all
nonzero polar Fourier coefficients through
\begin{equation}\label{eq:T-reconstruction-coefficients}
a_n(r)=\frac{y_{T,n}(r)}r,
\qquad
b_n(r)=\frac{r\delta_n(r)-y_{T,n}'(r)}{in},
\qquad n\ne0.
\end{equation}
The sign in the second formula also agrees with the
rotation identity
\[
\widetilde a_n
=\frac{y_{T,n}'-r\delta_n}{in}=-b_n.
\]

\subsection{The zeroth mode and the exact kernel}

Setting $n=0$ in \eqref{eq:div-mode} gives
$(ra_0)'(r)=r\delta_0(r).$
Since the field vanishes near the inner boundary, the
integration constant is zero. Hence
\begin{equation}\label{eq:T-a0}
a_0(r)=\frac1r\int_\eps^r s\delta_0(s)\,ds.
\end{equation}
The longitudinal reconstruction for $\vg$ leaves only
$\widetilde a_0=-b_0$ undetermined. Thus the field is
reconstructed in the form
\begin{equation}\label{eq:T-synthesis}
\vf(r,\theta)
=\sum_{n\in\Z}a_n(r)e^{in\theta}e_r(\theta)
+b_0(r)e_\theta(\theta)
+\sum_{n\ne0}b_n(r)e^{in\theta}e_\theta(\theta),
\end{equation}
where the transverse weighted pair determines every
coefficient except $b_0(r)$.
\begin{proof}[Proof of Theorem~\ref{thm:Tpair}]
Let $\vf\in\Vc_{\eps,R}$ and set $\vg=J\vf$.
Since $J^{-1}=-J$, the map $\vf\mapsto J\vf$ is a
bijection of $\Vc_{\eps,R}$ onto itself. By
\eqref{eq:transverse-rotation},
$\Tc\vf=\Lc\vg,
$ and $
\Tc^1\vf=\Lc^1\vg.$
Thus the transverse weighted data of $\vf$ determine the
longitudinal weighted data of $\vg$.

By Theorem~\ref{thm:Lpair}, these data determine
$\widetilde b_n$ for every $n\in\Z$ and
$\widetilde a_n$ for every $n\ne0$. The relations
$\widetilde a_n=-b_n,
\,
\widetilde b_n=a_n$
therefore determine $a_n$ for every $n\in\Z$ and $b_n$
for every $n\ne0$. Explicitly,
\[
a_n(r)=\frac{y_{T,n}(r)}r,
\qquad
b_n(r)=\frac{r\delta_n(r)-y_{T,n}'(r)}{in},
\qquad n\ne0,
\]
while
\[
a_0(r)=\frac1r\int_\eps^r s\delta_0(s)\,ds.
\]

To identify the kernel, suppose first that
$\Tc\vf=\Tc^1\vf=0$. Then
$\Lc\vg=\Lc^1\vg=0$, so Theorem~\ref{thm:Lpair} gives
$\vg=u(r)e_r(\theta)$ for some $u\in C_c^\infty((\eps,R)).$
Since $Je_r=e_\theta$, we obtain
$\vf=-J\vg=-u(r)e_\theta(\theta).
$
Hence
\[
\ker(\Tc,\Tc^1)
\subseteq
\{v(r)e_\theta(\theta):
v\in C_c^\infty((\eps,R))\}.
\]

Conversely, let $\vf=v(r)e_\theta(\theta)$ with
$v\in C_c^\infty((\eps,R))$. Since $Je_\theta=-e_r$,
\[
J\vf=-v(r)e_r(\theta)\in\ker(\Lc,\Lc^1)
\]
by Theorem~\ref{thm:Lpair}. Therefore,
\[
\Tc\vf=\Lc(J\vf)=0,
\qquad
\Tc^1\vf=\Lc^1(J\vf)=0.
\]
This proves the reverse inclusion and establishes
\[
\ker(\Tc,\Tc^1)
=
\{v(r)e_\theta(\theta):
v\in C_c^\infty((\eps,R))\}.
\]
For exact data, the recovered coefficients coincide with
the corresponding coefficients of the original smooth
compactly supported field. Their Fourier series converge
with all derivatives on compact sets. Consequently, the data determine $\vf$ up to the undetermined
term $b_0(r)e_\theta(\theta)$.
\end{proof}
\section{Complete recovery using complementary angular averages}
\label{sec:averages}

The longitudinal weighted pair determines every polar Fourier
coefficient except $a_0$, while the transverse weighted pair
determines every coefficient except $b_0$. We show that the
angular average of the complementary unweighted transform
determines the remaining coefficient.

\subsection{Rotation covariance and averaging}

Let
$O_\gamma=
\begin{pmatrix}
\cos\gamma&-\sin\gamma\\
\sin\gamma&\cos\gamma
\end{pmatrix}$
be the counterclockwise rotation matrix through angle
$\gamma$. Define its action on vector fields by
\begin{equation}\label{eq:rotation-action}
(U_\gamma\vf)(\vx)
:=O_{-\gamma}\vf(O_\gamma\vx).
\end{equation}
Thus $U_\gamma$ rotates the spatial argument by $\gamma$
and the vector value by $-\gamma$.

For a data function $G(\rho,\phi)$, denote its angular
average by
\[
\overline G(\rho)
:=\frac1{2\pi}\int_0^{2\pi}G(\rho,\phi)\,d\phi.
\]
The following lemma establishes the rotation covariance
of the circular transforms and shows that averaging
$U_\gamma\vf$ over $\gamma$ retains exactly the zeroth
Fourier coefficients of the polar components.

\begin{lemma}[Rotation covariance and angular averaging]
\label{lem:rotation-covariance}
For every $\vf\in\Vc_{\eps,R}$, $\gamma\in\R$, and
$(\rho,\phi)\in\Yc$, the circular transforms satisfy
\begin{align}
\Lc(U_\gamma\vf)(\rho,\phi)
&=\Lc\vf(\rho,\phi+\gamma),
\label{eq:rotation-L}\\
\Tc(U_\gamma\vf)(\rho,\phi)
&=\Tc\vf(\rho,\phi+\gamma),
\label{eq:rotation-T}
\end{align}
where the angular variable is understood modulo $2\pi$.

Moreover, the averaged field
\begin{equation}\label{eq:field-average}
\vf^{(0)}(\vx)
:=\frac1{2\pi}\int_0^{2\pi}
(U_\gamma\vf)(\vx)\,d\gamma
=a_0(r)e_r(\theta)+b_0(r)e_\theta(\theta),
\qquad \vx=r e_r(\theta),
\end{equation}
belongs to $\Vc_{\eps,R}$. Its longitudinal and transverse
transforms satisfy
\begin{equation}\label{eq:averaged-transforms}
\Lc\vf^{(0)}(\rho,\phi)=\overline{\Lc\vf}(\rho),
\qquad
\Tc\vf^{(0)}(\rho,\phi)=\overline{\Tc\vf}(\rho),
\qquad (\rho,\phi)\in\Yc.
\end{equation}
\end{lemma}

\begin{proof}
The change of variables $\vy=O_\gamma\vx$ maps
$C(\rho,\phi)$ onto $C(\rho,\phi+\gamma)$ and preserves
arclength. Indeed,
$O_\gamma\vc_\phi=\vc_{\phi+\gamma},
$ and $
|\vy-\vc_{\phi+\gamma}|
=|\vx-\vc_\phi|=\rho.$
The normal and tangent vectors transform according to
$O_\gamma\vn_\phi(\vx)
=\vn_{\phi+\gamma}(\vy),
$ and $
O_\gamma\vtau_\phi(\vx)
=\vtau_{\phi+\gamma}(\vy).$
Consequently,
\begin{align*}
\Lc(U_\gamma\vf)(\rho,\phi)
&=\int_{C(\rho,\phi)}
O_{-\gamma}\vf(O_\gamma\vx)
\cdot\vtau_\phi(\vx)\,ds_{\vx}\\
&=\int_{C(\rho,\phi)}
\vf(O_\gamma\vx)\cdot
O_\gamma\vtau_\phi(\vx)\,ds_{\vx}\\
&=\int_{C(\rho,\phi+\gamma)}
\vf(\vy)\cdot\vtau_{\phi+\gamma}(\vy)\,ds_{\vy}\\
&=\Lc\vf(\rho,\phi+\gamma).
\end{align*}
The same calculation with $\vn_\phi$ in place of
$\vtau_\phi$ proves \eqref{eq:rotation-T}.

For $\vx=r e_r(\theta)$, we have
$O_\gamma\vx=r e_r(\theta+\gamma)$,
$O_{-\gamma}e_r(\theta+\gamma)=e_r(\theta),
$ and $
O_{-\gamma}e_\theta(\theta+\gamma)=e_\theta(\theta).$
It follows that
\[
(U_\gamma\vf)(r,\theta)
=f_r(r,\theta+\gamma)e_r(\theta)
+f_\theta(r,\theta+\gamma)e_\theta(\theta).
\]
By periodicity, we have
\begin{align*}
\frac1{2\pi}\int_0^{2\pi}
f_r(r,\theta+\gamma)\,d\gamma=a_0(r),\qquad \text{and} \qquad
\frac1{2\pi}\int_0^{2\pi}
f_\theta(r,\theta+\gamma)\,d\gamma=b_0(r).
\end{align*}
This proves \eqref{eq:field-average}. The coefficients
$a_0,b_0$ are smooth and supported in a common compact
subinterval of $(\eps,R)$, so
$\vf^{(0)}\in\Vc_{\eps,R}$.

Finally, interchanging the integrals and using
\eqref{eq:rotation-L}, we obtain
\begin{align*}
\Lc\vf^{(0)}(\rho,\phi)
&=\frac1{2\pi}\int_0^{2\pi}
\Lc(U_\gamma\vf)(\rho,\phi)\,d\gamma\\
&=\frac1{2\pi}\int_0^{2\pi}
\Lc\vf(\rho,\phi+\gamma)\,d\gamma\\
&=\overline{\Lc\vf}(\rho).
\end{align*}
The corresponding argument for $\Tc$ proves the second
identity in \eqref{eq:averaged-transforms}.
\end{proof}

\subsection{Recovery of the undetermined coefficients}

By \eqref{eq:field-average},
\[
\vf^{(0)}(r,\theta)
=a_0(r)e_r(\theta)+b_0(r)e_\theta(\theta).
\]
Since $a_0$ and $b_0$ are independent of $\theta$, the
polar divergence and curl formulas give
\begin{align}
d_0(r)
&:=\diver\vf^{(0)}
=\frac1r(ra_0)'(r),
\label{eq:d0}\\
\omega_0(r)
&:=\curl\vf^{(0)}
=\frac1r(rb_0)'(r).
\label{eq:omega0-average}
\end{align}
Here $d_0$ and $\omega_0$ are regarded as radial scalar
fields on the plane, extended by zero outside the annulus.
Applying the fundamental scalar identities to $\vf^{(0)}$
and using \eqref{eq:averaged-transforms}, we obtain
\begin{align}
\partial_\rho\overline{\Tc\vf}(\rho)
&=\Rc d_0(\rho,\phi),
\label{eq:average-T-scalar}\\
\partial_\rho\overline{\Lc\vf}(\rho)
&=\Rc\omega_0(\rho,\phi).
\label{eq:average-L-scalar}
\end{align}
Because both scalar sources are radial, their transforms
are independent of $\phi$. Thus the zeroth Fourier--Abel
inverse determines $d_0$ from $\overline{\Tc\vf}$ and
$\omega_0$ from $\overline{\Lc\vf}$.

Multiplying \eqref{eq:d0} by $r$ and integrating gives
\[
ra_0(r)-\eps a_0(\eps)
=\int_\eps^r s d_0(s)\,ds.
\]
The zero extension of $a_0$ vanishes near $\eps$, so
$a_0(\eps)=0$. Applying the same argument to
\eqref{eq:omega0-average}, we obtain
\begin{align}
a_0(r)
&=\frac1r\int_\eps^r s d_0(s)\,ds,
\label{eq:average-recover-a0}\\
b_0(r)
&=\frac1r\int_\eps^r s\omega_0(s)\,ds.
\label{eq:average-recover-b0}
\end{align}

\begin{proof}[Proof of Theorem~\ref{thm:averages}]
Suppose first that
$\Lc\vf,$ $ \Lc^1\vf,$ and $
\overline{\Tc\vf}$
are known. Theorem~\ref{thm:Lpair} determines every polar
Fourier coefficient except $a_0$. Equation
\eqref{eq:average-T-scalar} determines $\Rc d_0$, and
Lemma~\ref{lem:scalar-inversion} reconstructs $d_0$.
Formula~\eqref{eq:average-recover-a0} then determines $a_0$.
Thus \eqref{eq:L-synthesis} reconstructs $\vf$ completely.

Similarly, suppose that
$\Tc\vf,$ $ \Tc^1\vf,$ and $
\overline{\Lc\vf}$
are known. Theorem~\ref{thm:Tpair} determines every
coefficient except $b_0$. Equation
\eqref{eq:average-L-scalar} and scalar inversion determine
$\omega_0$, after which
\eqref{eq:average-recover-b0} determines $b_0$.
The complete field follows from \eqref{eq:T-synthesis}.

For uniqueness in the first case, let $\vf_1,\vf_2$
produce the same data and set $\boldsymbol{h}:=\vf_1-\vf_2$.
Theorem~\ref{thm:Lpair} gives
\[
\boldsymbol{h}=u(r)e_r(\theta),
\qquad u\in C_c^\infty((\eps,R)).
\]
Since $\overline{\Tc\boldsymbol{h}}=0$, equation
\eqref{eq:average-T-scalar} and scalar injectivity imply
$\dfrac1r(ru(r))'=0.$
Hence $ru(r)$ is constant. Because $u$ vanishes near
$\eps$, this constant is zero, and therefore $\boldsymbol{h}=0$.

In the second case, Theorem~\ref{thm:Tpair} gives
\[
\boldsymbol{h}=v(r)e_\theta(\theta),
\qquad v\in C_c^\infty((\eps,R)).
\]
The equality $\overline{\Lc\boldsymbol{h}}=0$ similarly implies
$\dfrac1r(rv(r))'=0,$
and compact support gives $v=0$.

For exact data, the reconstruction recovers all Fourier
coefficients of the original smooth compactly supported
field. The Fourier syntheses therefore converge with all
derivatives on compact sets.
\end{proof}

\begin{remark}
Each supplementary measurement is a function of the
measurement radius $\rho$. A single additional scalar
linear measurement cannot eliminate the infinite-dimensional
kernel of either weighted pair. The theorem gives a
sufficient supplementary data set without asserting
minimality among all acquisition designs.
\end{remark}
\section{Recovery using two independent linear weights}
\label{sec:weights}

We prove Proposition~\ref{thm:two-weights} by reducing the
reconstruction of the vector field to scalar circular inversion
of its Cartesian components.

\subsection{The weighted disk identities}

Let $b\in C^\infty(S^1;\R^2)$ and write
$b(\phi)=\bigl(b^{(1)}(\phi),b^{(2)}(\phi)\bigr).$
For each fixed $\phi$, the weight
$w_b(\vx,\phi)=(\vx-\vc_\phi)\cdot b(\phi)$
satisfies
$\partial_1w_b=b^{(1)}(\phi)
$ and $
\partial_2w_b=b^{(2)}(\phi).$
Set $\omega=\curl\vf$. By the product rule,
\begin{align}
\curl(w_b\vf)
&=\partial_1(w_bf_2)-\partial_2(w_bf_1)\nonumber\\
&=w_b(\partial_1f_2-\partial_2f_1)
  +(\partial_1w_b)f_2-(\partial_2w_b)f_1\nonumber\\
&=w_b\omega+b^{(1)}(\phi)f_2-b^{(2)}(\phi)f_1\nonumber\\
&=w_b\omega+\vf\cdot Jb(\phi).
\label{eq:general-weight-product}
\end{align}
Applying Green's theorem to $w_b\vf$, we obtain
\begin{align}
\mathcal L_{C,b}\vf(\rho,\phi)
&=\int_{C(\rho,\phi)}
w_b(\vx,\phi)\vf(\vx)\cdot\vtau_\phi(\vx)\,ds_{\vx}
\nonumber\\
&=\int_{B(\rho,\phi)}w_b(\vx,\phi)\omega(\vx)\,d\vx
 +\int_{B(\rho,\phi)}\vf(\vx)\cdot Jb(\phi)\,d\vx.
\label{eq:general-weight-disk}
\end{align}

\subsection{Reconstruction of the Cartesian components}

\begin{proof}[Proof of Proposition~\ref{thm:two-weights}]
Suppose $b_1,b_2\in C^\infty(S^1;\R^2)$ satisfy
$D(\phi):=\det\bigl(b_1(\phi),b_2(\phi)\bigr)\ne0$ for every $\phi\in S^1.$
First, the identity
$\partial_\rho\Lc\vf=\Rc\omega$
and Lemma~\ref{lem:scalar-inversion} determine
$\omega=\curl\vf$.

For $j=1,2$, define the known residual
\[
Q_j(\rho,\phi)
:=\mathcal L_{C,b_j}\vf(\rho,\phi)
-\int_{B(\rho,\phi)}
w_{b_j}(\vx,\phi)\omega(\vx)\,d\vx.
\]
By \eqref{eq:general-weight-disk}, we have
\[
Q_j(\rho,\phi)
=\int_{B(\rho,\phi)}\vf(\vx)\cdot Jb_j(\phi)\,d\vx.
\]
Introduce
$F(\rho,\phi)
:=
\begin{pmatrix}
\Dc f_1(\rho,\phi)\\
\Dc f_2(\rho,\phi)
\end{pmatrix}.$
Since $Jb_j(\phi)$ is independent of $\vx$, these identities
give the linear system
\begin{equation}\label{eq:weight-system}
\begin{pmatrix}
Q_1(\rho,\phi)\\
Q_2(\rho,\phi)
\end{pmatrix}
=M(\phi)F(\rho,\phi),
\qquad
M(\phi)
:=
\begin{pmatrix}
(Jb_1(\phi))^T\\
(Jb_2(\phi))^T
\end{pmatrix}.
\end{equation}
Writing
$b_j=(b_j^{(1)},b_j^{(2)})$, we have
\[
M(\phi)=
\begin{pmatrix}
-b_1^{(2)}(\phi)&b_1^{(1)}(\phi)\\
-b_2^{(2)}(\phi)&b_2^{(1)}(\phi)
\end{pmatrix},
\]
and therefore
\[
\det M(\phi)
=b_1^{(1)}(\phi)b_2^{(2)}(\phi)
-b_1^{(2)}(\phi)b_2^{(1)}(\phi)
=D(\phi)\ne0.
\]
Thus the disk transforms of both components are determined by
\[
F(\rho,\phi)
=M(\phi)^{-1}
\begin{pmatrix}
Q_1(\rho,\phi)\\
Q_2(\rho,\phi)
\end{pmatrix}.
\]
Moreover, continuity and compactness give
$\min_{\phi\in S^1}|D(\phi)|>0,$
so $M^{-1}$ is smooth on $S^1$.

Since $M$ is independent of $\rho$, differentiating and applying
Lemma~\ref{lem:disk-differentiation} yields
\[
\begin{pmatrix}
\Rc f_1(\rho,\phi)\\
\Rc f_2(\rho,\phi)
\end{pmatrix}
=\partial_\rho F(\rho,\phi)
=M(\phi)^{-1}
\begin{pmatrix}
\partial_\rho Q_1(\rho,\phi)\\
\partial_\rho Q_2(\rho,\phi)
\end{pmatrix}.
\]
Both $f_1$ and $f_2$ belong to $C_c^\infty(A(\eps,R))$.
Lemma~\ref{lem:scalar-inversion} therefore reconstructs each
component uniquely from its scalar circular transform.
Hence the three data sets
$\Lc\vf,$ $
\mathcal L_{C,b_1}\vf$ and $
\mathcal L_{C,b_2}\vf$,
determine $\vf$ uniquely and constructively.
\end{proof}

\section{Examples and reconstruction procedures}
\label{sec:examples}

\subsection{Explicit smooth fields in the weighted kernels}

\begin{example}[A nonzero radial kernel field]
\label{ex:radial-kernel}
Choose $\eps<r_1<r_2<R$ and define
\[
u(r)=
\begin{cases}
\displaystyle
\exp\left(-\frac1{(r-r_1)(r_2-r)}\right),
& r_1<r<r_2,\\[5pt]
0,&\text{otherwise}.
\end{cases}
\]
The function $u$ is strictly positive on $(r_1,r_2)$, and
all its derivatives vanish at $r_1$ and $r_2$. Hence
$u\in C_c^\infty((\eps,R))$ and $u\ne0$.

Define
\[
\vf(\vx)=
\begin{cases}
\displaystyle u(|\vx|)\frac{\vx}{|\vx|},&\vx\ne0,\\[5pt]
0,&\vx=0.
\end{cases}
\]
Since $u$ vanishes near the origin and outside $[r_1,r_2]$,
we have
\[
\vf\in\Vc_{\eps,R},
\qquad
\supp\vf\subseteq
\{\vx:r_1\le|\vx|\le r_2\}.
\]
In polar coordinates, $\vf=u(r)e_r(\theta)$, so
\[
\curl\vf=0,
\qquad
\diver\vf=\frac1r(ru(r))'
=u'(r)+\frac{u(r)}r.
\]
By Theorem~\ref{thm:Lpair}, we have
$\Lc\vf=\Lc^1\vf=0.$
Thus, the longitudinal weighted pair cannot distinguish
$\vf$ from the zero field.

Nevertheless, the complementary averaged transverse data
determine
$d_0(r):=u'(r)+\dfrac{u(r)}r$
through \eqref{eq:average-T-scalar} and the zeroth
Fourier--Abel inverse. Since $u$ vanishes near $\eps$,
\[
\frac1r\int_\eps^r s d_0(s)\,ds
=\frac1r\int_\eps^r(su(s))'\,ds
=u(r).
\]
Hence the complementary averaged measurement recovers
the radial coefficient left undetermined by the
longitudinal weighted pair.
\end{example}

\begin{example}[A nonzero tangential kernel field]
\label{ex:tangential-kernel}
Using the same function $u$, define
\[
\vg(\vx)=
\begin{cases}
\displaystyle
u(|\vx|)\frac{(-x_2,x_1)}{|\vx|},&\vx\ne0,\\[5pt]
0,&\vx=0.
\end{cases}
\]
Then $\vg\in\Vc_{\eps,R}$ and
$\vg=u(r)e_\theta(\theta)$. The polar divergence and curl
formulas give
\[
\diver\vg=0,
\qquad
\curl\vg=\frac1r(ru(r))'
=u'(r)+\frac{u(r)}r.
\]
By Theorem~\ref{thm:Tpair}, we have
$\Tc\vg=\Tc^1\vg=0.$

The complementary averaged longitudinal data determine
$\omega_0(r):=u'(r)+\frac{u(r)}r$
through \eqref{eq:average-L-scalar} and the zeroth
Fourier--Abel inverse. Formula~\eqref{eq:average-recover-b0}
then gives
\[
u(r)=\frac1r\int_\eps^r s\omega_0(s)\,ds.
\]
Thus the complementary averaged measurement recovers
the tangential coefficient left undetermined by the
transverse weighted pair.

More generally, the maps
\[
u\longmapsto u(r)e_r(\theta),
\qquad
u\longmapsto u(r)e_\theta(\theta)
\]
are injective linear maps from $C_c^\infty((\eps,R))$
into the longitudinal and transverse weighted kernels,
respectively. Since $C_c^\infty((\eps,R))$ is infinite
dimensional, both weighted kernels are infinite dimensional.
\end{example}

\subsection{Exact-data reconstruction procedures}
The two weighted inverses use the same scalar routine. For data from an admissible field, perform the following operations, taking the longitudinal or transverse choice consistently:
\begin{enumerate}[label=(\arabic*)]
\item Differentiate the unweighted data in $\rho$ and
apply Lemma~\ref{lem:scalar-inversion} mode by mode.
This reconstructs $\omega=\curl\vf$ from $\Lc\vf$,
or $\delta=\diver\vf$ from $\Tc\vf$.
\item Subtract the weighted disk integral of that known source from the corresponding weighted data, forming \(Q_L\) in \eqref{eq:QL} or \(Q_T\) in \eqref{eq:QT}. Differentiate to obtain \(P_L\) or \(P_T\).
\item For \(n\ne0\), use \eqref{eq:GL} or \eqref{eq:GT} to compute the corrected coefficient. Apply the scalar inverse to \(2\rho G_{L,n}\) or \(2\rho G_{T,n}\), including this normalization factor. The result is \(y_{L,n}=rb_n\) or \(y_{T,n}=ra_n\), respectively. Recover the individual coefficients by \eqref{eq:L-reconstruction-coefficients} or \eqref{eq:T-reconstruction-coefficients}.
\item Recover the known zeroth profile by \eqref{eq:L-b0} or \eqref{eq:T-a0}. If the complementary average is supplied, its radius derivative and the zeroth scalar inverse give \(d_0\) or \(\omega_0\); use \eqref{eq:average-recover-a0} or \eqref{eq:average-recover-b0} to recover the remaining profile. Otherwise that profile remains the exact kernel ambiguity.
\item Synthesize \eqref{eq:L-synthesis} or \eqref{eq:T-synthesis}. The Cartesian components are
\[
f_1=f_r\cos\theta-f_\theta\sin\theta,
\qquad f_2=f_r\sin\theta+f_\theta\cos\theta.
\]
\end{enumerate}

\section{Numerical implementation and reconstruction results}
\label{sec:numerical}
In this section, we examine the numerical implementation of the reconstruction
formulas obtained in the preceding sections using MATLAB. We consider three sets of
measurements,
\begin{equation}\label{eq:num-data-pairs}
 (\Lc\vf,\Tc\vf),\qquad
 (\Lc\vf,\Lc^1\vf),\qquad
 (\Tc\vf,\Tc^1\vf).
\end{equation}
Our purpose is to examine finite-dimensional realizations of the exact-data
reconstruction formulas. The forward data are evaluated independently from
the prescribed vector field; inversion uses the second-kind Volterra equation
derived in \cite[Section~3.2]{AGL}. We present only noise-free
experiments to validate our inversion formulas numerically. For the two weighted pairs, the exact kernels identified in
Theorems~\ref{thm:Lpair} and~\ref{thm:Tpair} determine which representative
can be reconstructed.

\subsection{Phantoms and formation of the data}
\label{subsec:num-phantom}

We take $R=1$ and $\eps=0.2$. The field $\vf=(f_1,f_2)$ is defined on
$[-1,1]^2$ by
\begin{equation}\label{eq:num-phantom}
 f_j(\vx)=\chi(|\vx|)\sum_{k=1}^3 A_{jk}
 \exp\!\left(-\frac{|\vx-\boldsymbol\mu_{jk}|^2}{2\sigma_{jk}^2}\right),
 \qquad j=1,2,
\end{equation}
where the amplitudes, centers, and widths are given in
Table~\ref{tab:num-phantom}. The radial cutoff $\chi$ confines the phantom
to $0.27\leq |\vx|\leq0.89$; also consider
$[a,b]=[0.24,0.91]$. Here \(a\) and \(b\) denote the endpoints of the radial inversion interval, chosen so that the support of the phantom is contained in \((a,b)\).
\begin{table}[h!]
\centering
\begin{tabular}{@{}ccccc@{}}
\toprule
Component & $k$ & $A_{jk}$ & $\boldsymbol\mu_{jk}$ & $\sigma_{jk}$\\
\midrule
$f_1$ & 1 & $1.00$ & $(0.48,0.18)$ & $0.115$\\
$f_1$ & 2 & $-0.72$ & $(-0.40,-0.29)$ & $0.125$\\
$f_1$ & 3 & $0.48$ & $(-0.12,0.54)$ & $0.105$\\
$f_2$ & 1 & $0.95$ & $(-0.48,0.22)$ & $0.120$\\
$f_2$ & 2 & $-0.80$ & $(0.36,-0.38)$ & $0.125$\\
$f_2$ & 3 & $0.54$ & $(0.43,0.38)$ & $0.115$\\
\bottomrule
\end{tabular}
\caption{Parameters of the Gaussian phantoms.}
\end{table}\label{tab:num-phantom}
We use the same radial cutoff in all three configurations. Define
\begin{equation}\label{eq:num-smooth-cutoff}
 S(t)=\begin{cases}
 0,&t\leq 0,\\
 \displaystyle\frac{\exp(-1/t)}{\exp(-1/t)+\exp(-1/(1-t))},
 &0<t<1,\\
 1,&t\geq1,
 \end{cases}
 \qquad
 \chi(r)=S\!\left(\frac{r-0.27}{0.06}\right)
 S\!\left(\frac{0.89-r}{0.07}\right).
\end{equation}
Note that $\chi\in C_c^\infty((0.2,1))$ with
$\supp\chi\subset[0.27,0.89]$ and $\chi=1$ on
$[0.33,0.82]$. Consequently, the vector field in
\eqref{eq:num-phantom} belongs to $C_c^\infty(A(\eps,R);\R^2)$. 
The forward measurements are evaluated directly from the definitions of
the four circular transforms. For $j=0,\ldots,N_\phi-1$,
$k=0,\ldots,N_\rho-1$, and $\ell=0,\ldots,N_\alpha-1$, set
\[
 \phi_j=\frac{2\pi j}{N_\phi},\qquad
 \rho_k=\frac{k(R-a)}{N_\rho-1},\qquad
 \alpha_\ell=\frac{2\pi\ell}{N_\alpha},\qquad
 \vx_{kj\ell}=Re_{\phi_j}+\rho_k e_{\alpha_\ell}.
\]
Writing $\vn_\ell=e_{\alpha_\ell}$,
$\vtau_\ell=Je_{\alpha_\ell}$, and
$w_{kj\ell}=\rho_k\cos(\alpha_\ell-\phi_j)$, the periodic
trapezoidal rule gives
\begin{align}
 L_{kj}&=\frac{2\pi\rho_k}{N_\alpha}
  \sum_{\ell=0}^{N_\alpha-1}\vf(\vx_{kj\ell})\cdot\vtau_\ell,
 &T_{kj}&=\frac{2\pi\rho_k}{N_\alpha}
  \sum_{\ell=0}^{N_\alpha-1}\vf(\vx_{kj\ell})\cdot\vn_\ell,
 \label{eq:num-unweighted-data}\\
 L^1_{kj}&=\frac{2\pi\rho_k}{N_\alpha}
  \sum_{\ell=0}^{N_\alpha-1}w_{kj\ell}
  \vf(\vx_{kj\ell})\cdot\vtau_\ell,
 &T^1_{kj}&=\frac{2\pi\rho_k}{N_\alpha}
  \sum_{\ell=0}^{N_\alpha-1}w_{kj\ell}
  \vf(\vx_{kj\ell})\cdot\vn_\ell.
 \label{eq:num-weighted-data}
\end{align}
For evaluations, we used $N_r=321$ radial reconstruction nodes,
$N_\rho=601$ measurement radii, $N_\phi=160$ acquisition angles,
and $N_\alpha=1536$ quadrature points on each circle. For the scalar inversion, we employ 128-point Gauss–Legendre and 32-point Gauss–Chebyshev quadrature rules in all three reconstruction schemes. In the
results reported below, the angular Fourier series is truncated at
$|n|\leq8$ for all three data pairs. Since all of our reconstruction formula use the inversion of scalar circular Radon transform. So, for the sake of completeness, we first present the numerical implementation of the scalar inversion formula. 

\subsection{Numerical evaluation of the scalar Fourier--Abel inverse}
\label{subsec:num-volterra}
{We implement the scalar inversion recorded in
Section~\ref{subsec:scalar-inversion}, following
\cite[Section~3.2]{AGL}.}
The resulting second-kind Volterra equations are
discretized using piecewise-linear product integration
and solved by backward substitution.\\
Let $h$ be a scalar function with angular Fourier
coefficients $h_n$, and let $g=\Rc h$.
The Fourier coefficients $g_n(\rho_k)$ of the measured
data are approximated by the discrete Fourier transform
with respect to the acquisition angle:
\begin{equation}\label{eq:num-dft}
 g_n(\rho_k)\approx\frac1{N_\phi}\sum_{j=0}^{N_\phi-1}
 g(\rho_k,\phi_j)e^{-in\phi_j}.
\end{equation}
We reconstruct the nonnegative modes $0\leq n\leq N$; the negative
modes follow by complex conjugation for real-valued fields. Recall that
{with $t=R-\rho$, equation~\eqref{eq:abel-first}
takes the Abel form}
\begin{equation}\label{eq:num-abel}
 G_n(t):=g_n(R-t)=\int_t^b
 \frac{K_n(t,r)h_n(r)}{\sqrt{r-t}}\,dr.
\end{equation}
We use the kernel $K_n$ and diagonal normalization $d(r)=K_n(r,r)>0$
from \eqref{eq:abel-kernel} and \eqref{eq:abel-diagonal}.
{The Volterra equation
\eqref{eq:volterra-second-kind} and recovery formula
\eqref{eq:scalar-recovered} are implemented as}
\begin{equation}\label{eq:num-volterra}
 u_n(s)-\int_s^b V_n(s,r)u_n(r)\,dr=\gamma_n(s),
 \qquad h_n(s)=\frac{u_n(s)}{d(s)},
\end{equation}
where
\[
 V_n(s,r)=\frac1\pi\partial_sH_n(s,r),\qquad
 \gamma_n(s)=-\frac1\pi\frac{d}{ds}
 \int_s^b\frac{G_n(t)}{\sqrt{t-s}}\,dt.
\]
As the exact coefficients vanish near $b$, we have $G_n(b)=0$.
The expression for $\gamma_n$ is evaluated in the equivalent form
\begin{equation}\label{eq:num-gamma}
 \gamma_n(s)=-\frac1\pi\int_s^b
 \frac{G_n'(t)}{\sqrt{t-s}}\,dt
 =-\frac2\pi\int_0^{\sqrt{b-s}}G_n'(s+v^2)\,dv.
\end{equation}
Here a cubic spline represents the sampled function $G_n$, and
Gauss--Legendre quadrature evaluates the nonsingular integral on the
right. Substituting $
t=s+(r-s)v$, the Volterra kernel becomes
\begin{equation}\label{eq:num-kernel-evaluation}
 V_n(s,r)=\frac1{\pi d(r)}\int_0^1
 \frac{(1-v)\partial_tK_n(s+(r-s)v,r)}{\sqrt{v(1-v)}}\,dv,
\end{equation}
which is evaluated using Gauss--Chebyshev quadrature for the displayed weight. We discretize \eqref{eq:num-volterra} on
$a=s_0<s_1<\cdots<s_M=b$. On each interval $[s_j,s_{j+1}]$ we
approximate $u_n$ by the linear interpolant
\[
 u_n(s_j+\xi\Delta s_j)\approx(1-\xi)u_{n,j}+\xi u_{n,j+1},
 \qquad\Delta s_j=s_{j+1}-s_j.
\]
For the two-point Gauss nodes $\xi_\pm=(1\pm1/\sqrt3)/2$, define
\begin{align}
 A_{ij}&=\frac{\Delta s_j}{2}\sum_{q\in\{-,+\}}
  (1-\xi_q)V_n(s_i,s_j+\xi_q\Delta s_j),\nonumber\\
 B_{ij}&=\frac{\Delta s_j}{2}\sum_{q\in\{-,+\}}
  \xi_q V_n(s_i,s_j+\xi_q\Delta s_j).
 \label{eq:num-product-integration}
\end{align}
The resulting discrete equation is
\begin{equation}\label{eq:num-backward-solve}
 (1-A_{ii})u_{n,i}=\gamma_n(s_i)
 +\sum_{j=i+1}^{M-1}A_{ij}u_{n,j}
 +\sum_{j=i}^{M-1}B_{ij}u_{n,j+1},
\end{equation}
with $u_{n,M}=0$. The unknowns are determined successively for
$i=M-1,M-2,\ldots,0$. We then set $h_n(s_i)=u_{n,i}/d(s_i)$ and use
angular Fourier synthesis.
\subsection{Reconstruction from the joint unweighted data}
\label{subsec:num-unweighted}
\begin{figure}[h!]
\centering
\includegraphics[width=\textwidth]{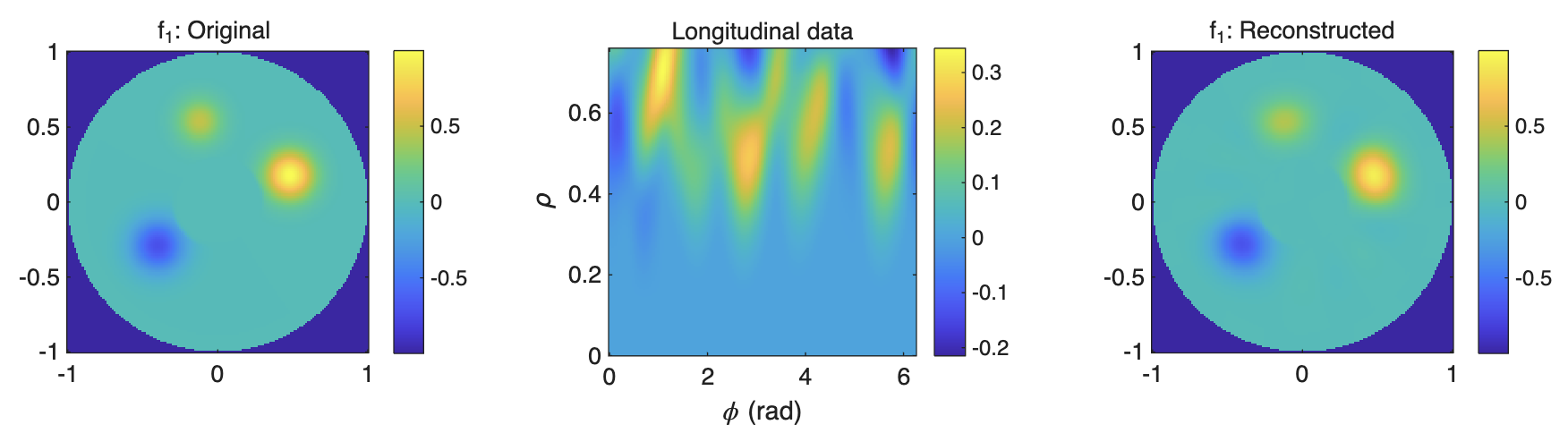}
\par\medskip
\includegraphics[width=\textwidth]{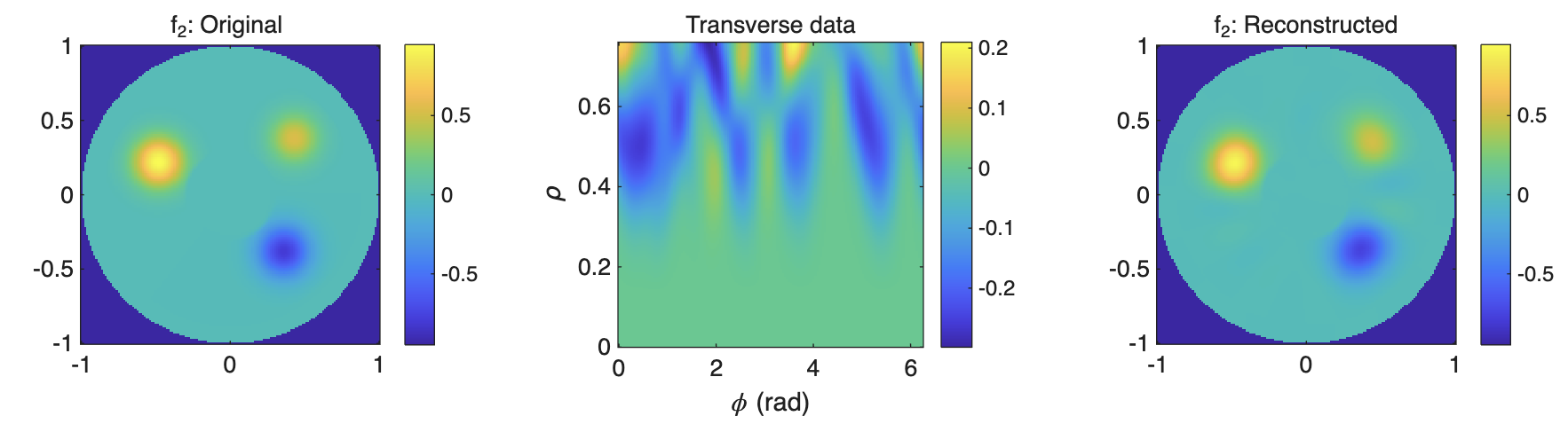}
\caption{Reconstruction of $\vf$ from $(\Lc\vf,\Tc\vf)$.
}
\label{fig:num-LT}
\end{figure}

The goal of this subsection is to recover $\vf$ from $(\Lc\vf,\Tc\vf)$, numerically and validate 
Theorem~\ref{thm:joint}. By \eqref{eq:fundamental-L} and
\eqref{eq:fundamental-T},
\[
 \partial_\rho\Lc\vf=\Rc\omega,\qquad
 \partial_\rho\Tc\vf=\Rc\delta,
 \qquad \omega=\curl\vf,\quad\delta=\diver\vf.
\]
After differentiating the angular Fourier coefficients of the two data
arrays in $\rho$, we apply \eqref{eq:num-volterra} to reconstruct the
coefficients of $\omega$ and $\delta$. Their truncated Fourier series
provide approximations of these scalar fields.
The components of $\vf$ are then recovered from the equations
\begin{equation}\label{eq:num-poisson}
 \Delta f_1=\partial_1\delta-\partial_2\omega,\qquad
 \Delta f_2=\partial_2\delta+\partial_1\omega.
\end{equation}
We approximate the first derivatives on a uniform grid in $[-1,1]^2$
by centered differences and solve the two Dirichlet problems by the
five-point discrete Laplacian, imposing zero values on the boundary of
the square. The reconstruction is performed using
$161\times161$ grid.
Figure~\ref{fig:num-LT} shows the phantom components, the
longitudinal and transverse data, and their reconstructions. 

\subsection{Reconstruction from the longitudinal and weighted longitudinal data}
\label{subsec:num-weighted-L}
The focus of this subsection is to recover $\vf$ from the data pair $(\Lc\vf,\Lc^1\vf)$. First using the scalar
inversion we recovered $\omega=\curl\vf$. From
\eqref{eq:PL} and \eqref{eq:GL}, we have $P_L=\partial_\rho\Lc^1\vf-\Rc(w_\phi\omega)$
and
\begin{equation}\label{eq:num-second-L}
 2\rho G_{L,n}(\rho)=RP_{L,n}(\rho)
  +\frac{2\rho}{n}\int_{R-\rho}^b
  r\omega_n(r)\sin(n\beta(r,\rho))\,dr,  \quad  n\ne0.
\end{equation}
To evaluate \(\mathcal R_C(w_\phi\omega)\), we use the identity
\[
w_\phi(x)=\frac{r^2-R^2-\rho^2}{2R},
\qquad r=|x|,\quad x\in C(\rho,\phi).
\]

The scalar circular integral is evaluated using the
substitution $r=R-\rho+v^2$, which removes the
Abel-type square-root singularity. The same
substitution is applied to the integral in
\eqref{eq:num-second-L}. Both resulting integrals
are nonsingular and are evaluated using
Gauss--Legendre quadrature.

\begin{figure}[h!]
\centering
\includegraphics[width=\textwidth]{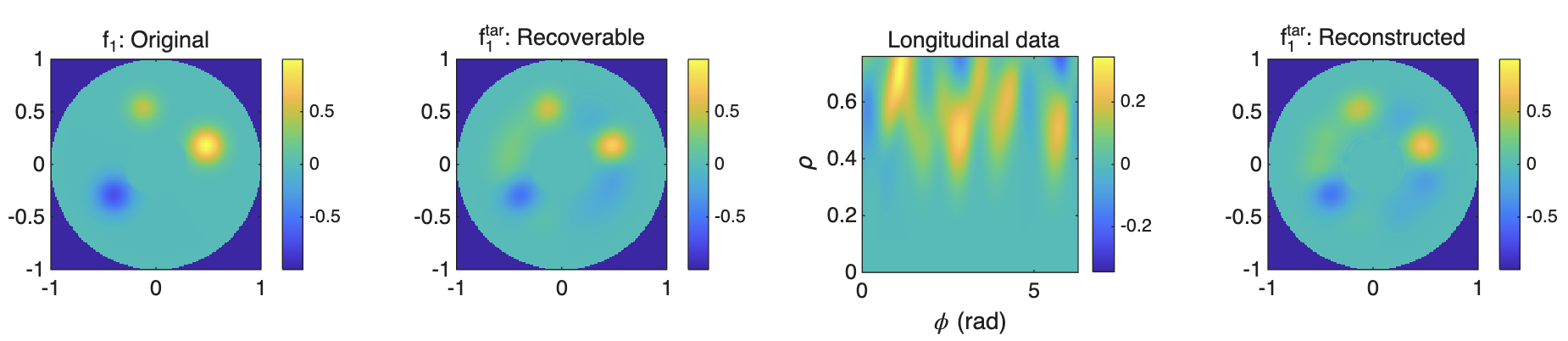}
\par\medskip
\includegraphics[width=\textwidth]{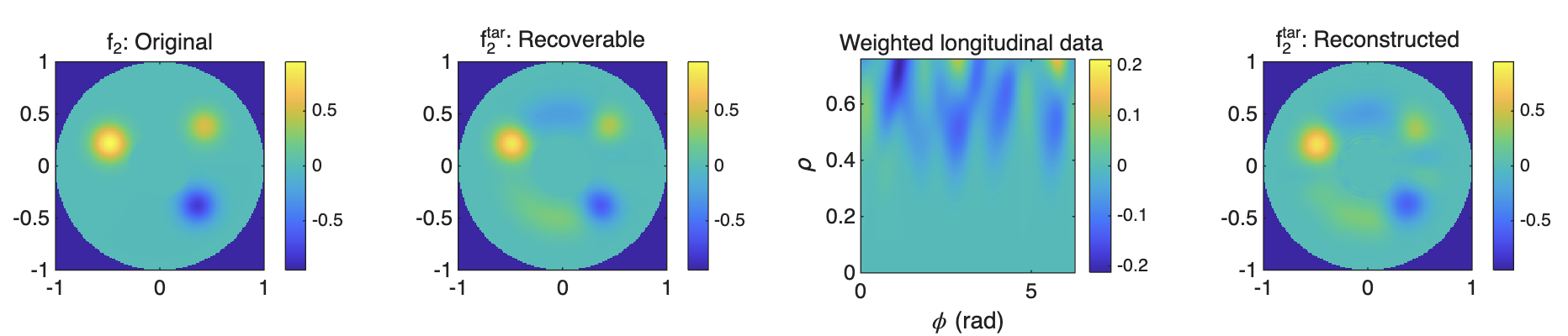}
\caption{Reconstruction of $\vf$ from $(\Lc\vf,\Lc^1\vf)$.
}
\label{fig:num-LL1}
\end{figure}

By \eqref{eq:L-scalar-normalization}, the quantity
$2\rho G_{L,n}$ is precisely the scalar circular-transform
coefficient of $y_{L,n}(r)e^{in\theta}$, where
$y_{L,n}=rb_n$. Thus the second application of the same Volterra
solver recovers $y_{L,n}$, after which
\begin{equation}\label{eq:num-L-coeff}
 b_n(r)=\frac{y_{L,n}(r)}r,\qquad
 a_n(r)=\frac{y_{L,n}'(r)-r\omega_n(r)}{in},\quad n\ne0,
 \qquad b_0(r)=\frac1r\int_\eps^r s\omega_0(s)\,ds.
\end{equation}
To compute $y'_{L,n}$, we fit a cubic polynomial
by least squares to the reconstructed values of
$y_{L,n}$ in a local radial window of half-width
$0.032$. The derivative of the fitted polynomial
is then evaluated at the corresponding grid point.
This reduces the amplification of numerical
oscillations during differentiation.
 Theorem~\ref{thm:Lpair} shows that the data cannot recover the
rotationally invariant radial component $a_0(r)e_r$. We choose
$a_0=0$ in the reconstructed representative and compare it with
\begin{equation}\label{eq:num-L-target}
 \vf_L^{\mathrm{tar}}(r,\theta)
 =\vf(r,\theta)-a_0(r)e_r(\theta),\qquad
 a_0(r)=\frac1{2\pi}\int_0^{2\pi}\vf(r,\theta)\cdot e_r(\theta)\,d\theta.
\end{equation}
The numerical reconstruction of the recoverable
field are presented in Figure~\ref{fig:num-LL1}.
\subsection{Reconstruction from the transverse and weighted transverse data}
\label{subsec:num-weighted-T} 
The objective of this subsection is to recover a vector field $\vf$ from the data set $(\Tc\vf,\Tc^1\vf)$.
The reconstruction from $(\Tc\vf,\Tc^1\vf)$ is performed using the
fixed rotation $\vg=J\vf$. Since
\[
 \Lc\vg=\Tc\vf,\qquad\Lc^1\vg=\Tc^1\vf,
 \qquad\curl\vg=\diver\vf=:\delta,
\]
the longitudinal procedure applies to $\vg$ without changing the
scalar inversion algorithm. 
First, we recover $\delta=\operatorname{div}\vf$ from
$\partial_\rho\Tc\vf=\Rc\delta$.
Using the recovered divergence and the weighted transverse
measurements, we compute $P_T=\partial_\rho\Tc^1\vf-\Rc(w_\phi\delta).$
The scalar data for the second inversion are then given by
\begin{equation}\label{eq:num-second-T}
2\rho G_{T,n}(\rho)
=RP_{T,n}(\rho)
+\frac{2\rho}{n}\int_{R-\rho}^b
r\delta_n(r)\sin(n\beta(r,\rho))\,dr,
\qquad n\ne0.
\end{equation}

The second Volterra inversion gives $y_{T,n}=ra_n$, and the polar
coefficients follow from
\begin{equation}\label{eq:num-T-coeff}
 a_n(r)=\frac{y_{T,n}(r)}r,\qquad
 b_n(r)=\frac{r\delta_n(r)-y_{T,n}'(r)}{in},\quad n\ne0,
 \qquad a_0(r)=\frac1r\int_\eps^r s\delta_0(s)\,ds.
\end{equation}
The same local-polynomial differentiation is used to compute
$y_{T,n}'$. Finally, we rotate the recovered field back by $\vf=-J\vg$. Here the undetermined term is $b_0(r)e_\theta$, by
Theorem~\ref{thm:Tpair}. The error reference is therefore
\begin{equation}\label{eq:num-T-target}
 \vf_T^{\mathrm{tar}}(r,\theta)
 =\vf(r,\theta)-b_0(r)e_\theta(\theta),\qquad
 b_0(r)=\frac1{2\pi}\int_0^{2\pi}
 \vf(r,\theta)\cdot e_\theta(\theta)\,d\theta.
\end{equation}

\begin{figure}[h!]
\centering
\includegraphics[width=\textwidth]{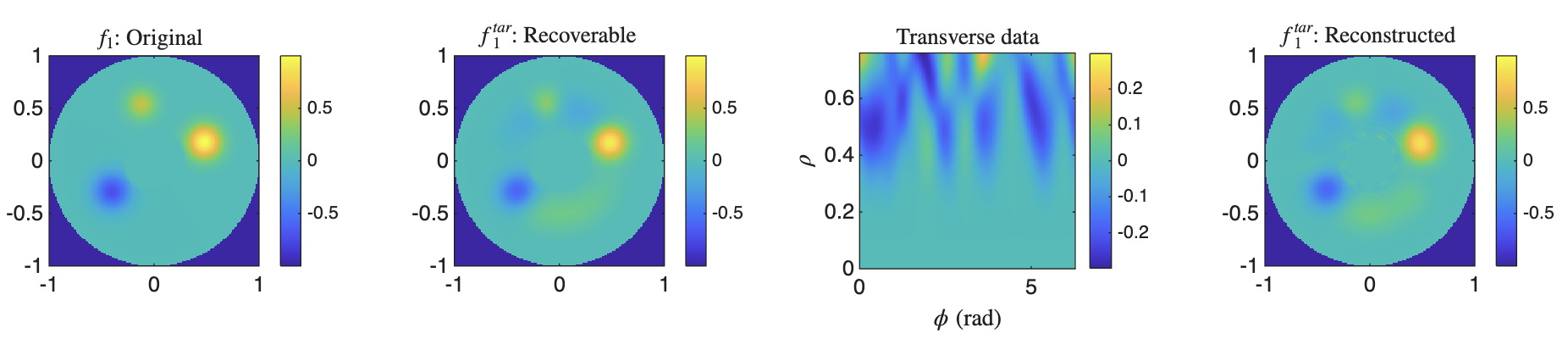}
\par\medskip
\includegraphics[width=\textwidth]{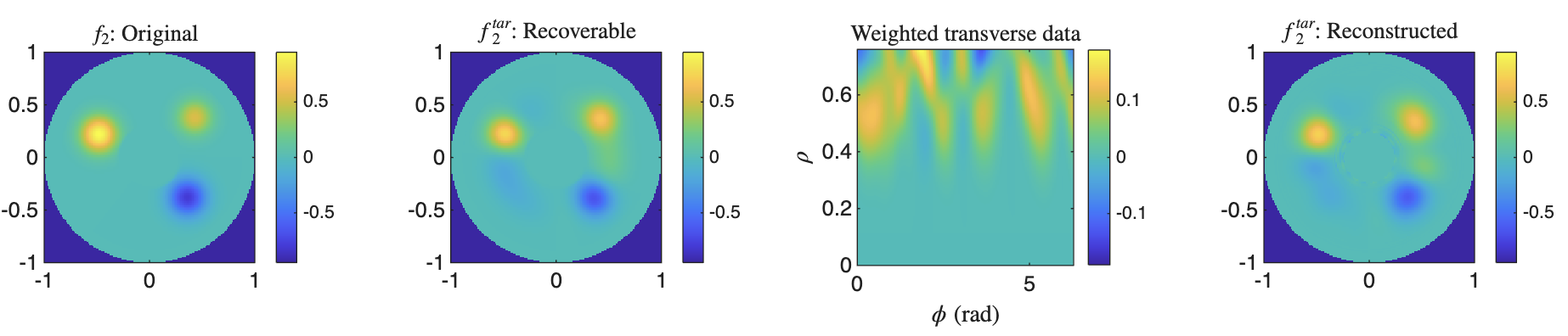}
\caption{Reconstruction of $\vf$ from $(\Tc\vf,\Tc^1\vf)$.}
\label{fig:num-TT1}
\end{figure}

The numerical reconstruction of the recoverable
field are presented in Figure~\ref{fig:num-TT1}. 
\subsection{Reconstruction errors}
\label{subsec:num-results}
We measure the accuracy of each recovered component. The relative $L^2$-error of the reconstructions,
calculated as follows:
\begin{equation}\label{eq:num-error}
 E_{f_j}=\frac{\|f_j^{\mathrm{rec}}-f_j^{\mathrm{tar}}\|_{2}}
 {\|f_j^{\mathrm{tar}}\|_{2}}\times100\%,\qquad j=1,2.
\end{equation}
where $\| h \|_2$ denotes the spectral norm, i.e., the maximum singular value of the matrix $h$ . The target is $\vf$ for the joint unweighted data, $\vf_L^{\mathrm{tar}}$
for the longitudinal weighted pair, and $\vf_T^{\mathrm{tar}}$ for
the transverse weighted pair. 
\begin{table}[h!]
\centering
\label{tab:num-errors}
\begin{tabular}{@{}lcc@{}}
\toprule
Measurements & $E_{f_1}$ & $E_{f_2}$\\
\midrule
$(\Lc\vf,\Tc\vf)$ & $8.3975\%$ & $9.3817\%$\\
$(\Lc\vf,\Lc^1\vf)$ & $11.0435\%$ & $11.6473\%$\\
$(\Tc\vf,\Tc^1\vf)$ & $10.8797\%$ & $11.4391\%$\\
\bottomrule
\end{tabular}
\caption{Componentwise relative reconstruction errors.}
\end{table}
\section*{Acknowledgements}\label{sec:acknowledgements}
R.\,B.\ gratefully acknowledges the financial support
provided by the Institute Postdoctoral Fellowship at the
Indian Institute of Technology Bombay (IIT Bombay).
 I.\,Z.\ acknowledge the financial support from RICAM, Austrian Academy of Sciences.
The first author also thanks Prof.\ Venky Krishnan for insightful
discussions and valuable suggestions concerning this work
during the ``Inverse Problems Students Meet'' held at
IISER Bhopal. The second author thanks Prof. Otmar Scherzer for valuable discussions regarding potential applications.

\vspace{.4cm}
 \noindent \textbf{Data availability statement.} \ No datasets were generated or analyzed during the current study; therefore, data sharing is not applicable.

    \vspace{.1cm}
\noindent \textbf{Conflict of interest.} \
The authors declare that they have no conflicts of interest regarding the research, authorship, and/or publication of this manuscript.
\enlargethispage{3\baselineskip}
\renewcommand{\baselinestretch}{1}\selectfont
\bibliographystyle{plain}
\bibliography{references}
\end{document}